\documentclass[11pt]{article}
\usepackage{amsmath}
\usepackage{dsfont} 
\usepackage{amssymb} 
\usepackage{graphicx,color}
\usepackage{extarrows}

\usepackage{mathrsfs} 

\begin{document}
\title{{\bf\Large Slow foliation of a slow-fast stochastic evolutionary system
\thanks{Part of this work was done while Guanggan Chen was participating
the program on ``Interactions Between Analysis and Geometry" at the
Institute for Pure \& Applied Mathematics, University of California,
Los Angeles, USA. This work was supported by the NSFC Grants
(Nos.11371267, 11071177, 11371367 and 11271290), the NSF Grant
1025422, and the Scientific Research Fund of Science and Technology
Bureau of Sichuan Province (Grant Nos. 2012JQ0041 and 2010JY0057). }
} }
\author{Guanggan Chen\\
College of Mathematics and Software Science,\\
Sichuan Normal University, Chengdu, 610068, China\\
\emph{E-mail: chenguanggan@hotmail.com}\\
\\
Jinqiao Duan\\
Institute for Pure and Applied Mathematics,\\
University of California,
Los Angeles, CA 90095, USA  \\
\& \\
Department of Applied Mathematics, \\
Illinois Institute of Technology,
Chicago, IL 60616, USA \\
\emph{E-mail: duan@iit.edu }\\
\\
Jian Zhang\\
College of Mathematics and Software Science, \\
Sichuan Normal University,
Chengdu, 610068, China\\
\emph{E-mail: zhangjiancdv@sina.com}}

\date{\today}
\maketitle


\begin{center}
\emph{Dedicated to Professor Zhiming Ma on the Occasion of his 65th
Birthday}
\end{center}

\baselineskip 7.5mm
\begin{center}
\begin{minipage}{120mm}
{\bf \large Abstract:} {\small This work is concerned with the
dynamics of a slow-fast stochastic evolutionary system quantified with a scale   parameter.
An invariant foliation decomposes the state space into geometric regions of different dynamical regimes,
and thus helps understand dynamics.  A slow invariant
foliation is established for this system. It is shown that the slow foliation converges to
a critical foliation (i.e., the scale parameter is zero) in probability distribution, as the
scale parameter tends to zero. The approximation of slow foliation is also
constructed with error estimate in distribution. Furthermore, the geometric
structure of the slow foliation is investigated: every fiber of
the slow foliation parallels each other, with the slow manifold as
a special fiber. In fact, when an arbitrarily chosen
  point of a fiber falls in the slow manifold, the fiber
  must be the slow manifold itself.

}


\end{minipage}
\end{center}


\renewcommand{\theequation}{\thesection.\arabic{equation}}
\setcounter{equation}{0}

\section{Introduction }

\quad\,   Random fluctuations may have delicate effects on dynamical evolution
of complex systems (\cite{Arnold, CDZ2, DZ, WD}). The slow-fast
stochastic evolutionary systems are appropriate mathematical models
for various multi-scale   systems under random  influences.
\par
We consider the following slow-fast stochastic evolutionary
system
\begin{equation}\label{Equation-s}
\frac{dx^\varepsilon}{dt}=Ax^\varepsilon+f(x^\varepsilon,
y^\varepsilon)+\sigma_1\dot{W_1},~~~~~~\quad \hbox{in}\quad H_s,
\end{equation}
\begin{equation}\label{Equation-f}
\frac{dy^\varepsilon}{dt}=\frac{1}{\varepsilon}By^\varepsilon+\frac{1}{\varepsilon}g(x^\varepsilon,
y^\varepsilon)+\frac{\sigma_2}{\sqrt{\varepsilon}}\dot{W_2},\quad
\hbox{in}\quad H_f,
\end{equation}
where   $\varepsilon$ is a small positive parameter
  ($0<\varepsilon\ll 1$).  The Hilbert spaces $H_s$ and $ H_f$,
linear operators $A$ and $ B$, nonlinearities $f$ and $g$, and mutually
independent Wiener processes $W_1$ and $W_2$ will be specified in
the next section. The white noises $\dot{W}_1$ and $\dot{W}_2$ are
the generalized time derivatives of $W_1$ and $W_2$, respectively.
The positive
constants $\sigma_1$ and $\sigma_2$ are the intensities of white
noises. Since the small  scale  parameter
$\varepsilon$   is such that $\|\frac{dx}{dt}\|_{H_s}\ll
\|\frac{dy}{dt}\|_{H_f}$, we usually say that $x$ is the ``slow" component and
$y$ is the ``fast" component.
\par

The main goal of this paper is to investigate state space decomposition for this system, by considering a
slow invariant foliation, and examining its approximation and structure.
\par

Invariant foliations and invariant manifolds play a
  significant  role in the study of the qualitative dynamical
behaviors, as they provide  geometric structures to understand or
reduce stochastic dynamics (\cite{BW, CDLS, CDZ, DLS-2003, DLS-2004,
Fenichel, LS, MZZ, RDJ}). An invariant foliation is about
quantifying certain sets (called fibers or leaves) in state space
for a dynamical system. A fiber consists of all those points
starting from which the dynamical orbits are exponentially
approaching each other, in forward time (``stable foliation") or
backward time (``unstable foliation"). These fibers are thus
building blocks for understanding dynamics, as they carry dynamical
information.  Collectively they provide a decomposition of the state
space.
\par

For a system  like  (\ref{Equation-s})-(\ref{Equation-f}), Schmalfuss and Schneider \cite{SS}
studied the slow manifold in the finite dimensional case. Wang,
Duan, and Roberts \cite{WRD-2012-JDE, WR-2013-JMAA} further studied the slow
manifold, and a relation with averaging as quantified via large deviations and approximations. In the infinite dimensional setting, Fu, Liu and
Duan \cite{FLD} investigated the slow manifold and its approximation. These research works are at
the level of geometric and global invariant sets. In the context of analyzing individual sample solution paths, Freidlin
\cite{Freidlin} used large deviation theory to describe the
dynamics, and Berglund and Gentz \cite{BG-book} showed that the sample solution
paths   are concentrated in a neighborhood of
the critical manifold (also see \cite{KP-book}).
\par

\par

Although invariant foliation theory has been developed for  deterministic systems
in \cite{BLZ, CHT, CLL},  it is still in   infancy for
stochastic evolutionary systems.
Recently, Lu and Schmalfuss \cite{LS08} studied the existence of
random invariant foliation for a class of stochastic partial differential equations, and Sun, Kan and Duan \cite{SKD} established
the approximation of random invariant foliations.


We define that a \emph{slow} foliation of a slow-fast system to be the foliation in
which the fibers are parameterized or represented  by slow variables, when the  scale parameter $\varepsilon$ is sufficiently small.
In a sense, the fast variables are eliminated.  A
\emph{critical} foliation corresponds to the foliation
  with   zero scale parameter.  Furthermore, the slow foliation converges to the critical
foliation, as the singular perturbation parameter $\varepsilon$
tends to zero.
\par


For system (\ref{Equation-s})-(\ref{Equation-f}), we establish the
existence of slow foliation, which is a graph of a Lipsichtz
continuous map. The dynamical orbits of
the slow-fast stochastic system are exponentially approaching each other in
backward time only if they start from the same slow fiber. In addition, we
show the slow foliation converges to a critical foliation in probability
distribution, as $\varepsilon$ tends to zero.
Furthermore, we examine the geometric structure of the slow foliation and show that   fibers of the slow foliation parallel with each
other. In fact, the slow manifold is one fiber of the slow foliation.
When an arbitrarily chosen point of the slow foliation is
in the slow manifold, the fiber   passing through the point is just the
slow manifold.


\par
This paper is organized as follows. In the next section, we make   hypotheses  for the slow-fast system and  recall
  basic concepts in random dynamical systems, including    random slow manifolds. In \S 3, we  present a motivating example about
slow foliation. In \S 4, we prove the existence of
slow foliation (Theorem 4.1), examine the geometric structure of the
slow foliation, and analyze a relationship between the slow
foliation and the slow manifold (Theorem 4.2). In \S 5,
we establish the existence of a critical foliation (Theorem 5.1), prove the convergence of the slow
foliation   to the critical foliation in probability distribution as the
scale parameter tends to zero (Theorem 5.2), and
construct an approximation of slow foliation in probability distribution
(Theorem 5.3).
\par

\renewcommand{\theequation}{\thesection.\arabic{equation}}
\setcounter{equation}{0}

\section{Preliminaries }

\subsection{Basic setup}
\quad\;   For the slow-fast system (\ref{Equation-s})-(\ref{Equation-f}), let
$H_s$ and $H_f$ be two separable Hilbert spaces with the
norms $\|\cdot\|_{H_s}$ and $\|\cdot\|_{H_f}$, respectively. The
space $H_s$ denotes the state space for slow variables, and $H_f$
the state space for fast variables. Henceforth,   we use
the subscripts or superscripts ``$s$" and ``$f$" to denote those spaces or quantities  that are
related to the slow variables and   fast variables, respectively.
We introduce the following hypotheses.
\par
{\bf Hypothesis H1 }(Dichotomy condition): The linear operator $A$
generates a $C_0$-semigroup $e^{At}$ on $H_s$ satisfying
$$
\|e^{At}x\|_{H_s}\leq e^{-\gamma_s t}\|x\|_{H_s},\quad
\hbox{for}\quad t\leq 0,
$$
and the linear operator $B$ generates a $C_0$-semigroup $e^{Bt}$ on $H_f$
satisfying
$$
\|e^{Bt}y\|_{H_f}\leq e^{-\gamma_f t}\|y\|_{H_f},\quad
\hbox{for}\quad t\geq 0,
$$
where $\gamma_s<0<\gamma_f$.
\par
{\bf Hypothesis H2} (Lipschitz condition): The nonlinear functions
$$
\begin{array}{l}
f: H_s\times H_f \longrightarrow H_s,\\
g: H_s\times H_f \longrightarrow H_f,
\end{array}
$$
  are $C^1$-smooth  with $f(0,0)=0$ and $g(0,0)=0$,
and  satisfy a   Lipschitz condition, i.e., there exists a positive
constant $K$ such that for every $(x,y)^T \in H_s\times H_f$ and
every $(\tilde{x},\tilde{y})^T \in H_s\times H_f$,
$$
\|f(x, y)-f(\tilde{x}, \tilde{y})\|_{H_s}\leq
K(\|x-\tilde{x}\|_{H_s}+\|y-\tilde{y}\|_{H_f}),
$$
$$
\|g(x, y)-g(\tilde{x}, \tilde{y})\|_{H_s}\leq
K(\|x-\tilde{x}\|_{H_s}+\|y-\tilde{y}\|_{H_f}).
$$
Here and hereafter, the superscript ``$T$" denotes the matrix
transpose.
\par

{\bf Hypothesis  H3 } (Gap condition): The Lipschitz constant $K$ of
the nonlinear functions $f$ and $g$ satisfies the condition
$K<\frac{-\gamma_s\cdot\gamma_f}{2\gamma_f-\gamma_s}$.
\par

\subsection{Random dynamical systems}

\quad \;
We recall some basic concepts in random dynamical
systems (\cite{DLS-2004}).
Let $(\Omega,\mathcal{F}, \mathds{P})$ be a probability space. A
flow $\theta$ of mappings $\{\theta_t\}_{t\in \mathbb{R}}$ is
defined on the sample space $\Omega$ such that
\begin{equation}
\theta: \mathbb{R}\times \Omega\to \Omega, \quad \theta_0=id,\quad
\theta_{t_1}\theta_{t_2}=\theta_{t_1+t_2},
\end{equation}
for $t_1, t_2 \in \mathbb{R}$. This flow is assumed to be
$(\mathcal{B}(\mathbb{R})\otimes\mathcal{F},
\mathcal{F})$-measurable, where $\mathcal{B}(\mathbb{R})$ is the
$\sigma$-algebra of Borel sets on the real line $\mathbb{R}$. To
have this measurability, it is not allowed to replace $\mathcal{F}$
by its $\mathds{P}$-completion $\mathcal{F}^{\mathds{P}}$; see
Arnold \cite[P547]{Arnold}. In addition, the measure $\mathds{P}$ is
assumed to be ergodic with respect to $\{\theta_t\}_{t\in
\mathbb{R}}$. Then $\Theta=(\Omega,\mathcal{F}, \mathds{P}, \theta)$
is called a \emph{metric dynamical system}.
\par
For our purpose, we will consider a special but very important
metric dynamical system induced by the Wiener process. Let $W(t)$ be
a two-sided Wiener process taking values in a Hilbert space $H$. Its sample paths are in the space
$C_0(\mathbb{R}, H)$ of  real continuous functions defined on
$\mathbb{R}$, taking zero value at $t = 0$. This set is equipped
with the compact open topology. On this set we consider the
measurable flow $\theta = \{\theta_t\}_{t\in \mathbb{R}}$, defined
by
$$\theta_t\omega =\omega(\cdot+t)-\omega(t),\quad \omega\in\Omega,\quad t\in \mathbb{R}.$$
The distribution
of this process induces a probability measure on
$\mathcal{B}(C_0(\mathbb{R}, H))$ and it is called the Wiener
measure. Note that this measure is ergodic with respect to  $\theta_t$; see \cite[Appendix A]{Arnold}. We also
consider, instead of the whole $C_0(\mathbb{R},H)$, a
$\{\theta_t\}_{t\in \mathbb{R}}$-invariant subset $\Omega\subset
C_0(\mathbb{R}, H)$) of $\mathds{P}$-measure one and the
trace $\sigma$-algebra $\mathcal{F}$ of
$\mathcal{B}(C_0(\mathbb{R},H))$ with respect to $\Omega$.
Recall that a set $\Omega$ is called {$\{\theta_t\}_{t\in \mathbb{R}}$-invariant
if $\theta_t\Omega = \Omega$ for $t \in \mathbb{R}$. On
$\mathcal{F}$, we consider the restriction of the Wiener measure
and still denote it by $\mathds{P}$.
\par
In general, the dynamics of a stochastic system on the state space
$H$ (often a Hilbert space) over the flow $\theta$ is described by a
cocycle.   A cocycle $\phi$ is a mapping:
$$
\phi:    \mathbb{R}\times\Omega\times H \to H,
$$
which is $(\mathcal{B}(\mathbb{R})\otimes \mathcal{F} \otimes
\mathcal{B}(H), \mathcal{F})$-measurable such that
$$
\begin{array}{l}
\phi(0, \omega, x) = x,\\
\phi(t_1 + t_2, \omega, x) = \phi(t_2, \theta_{t_1}\omega, \phi(t_1,
\omega, x)),
\end{array}
$$
for $t_1, t_2 \in \mathbb{R}, \omega \in \Omega$ and $x\in H$. Then
$\phi$ together with the metric dynamical system $\theta$ forms a
\emph{random dynamical system}.
\par
\emph{A stable fiber} and \emph{an unstable fiber} of a foliation
are defined as follows (also see \cite{CLL}).
\par
\par
(i)  \emph{$\mathcal{W}_{\beta s}(x,\omega)$ is called a
$\beta$-stable fiber} passing through $x\in H$ with $\beta\in
\mathbb{R}^-$,  if
$\|\phi(t,\omega,x)-\phi(t,\omega,\widetilde{x})\|_H=O(e^{\beta t}),
\forall\; \omega\in \Omega$ as $t\to +\infty$ for all  $x,
\widetilde{x}\in \mathcal{W}_{\beta s}$.
\par
(ii)   \emph{$\mathcal{W}_{\beta u}(x,\omega)$ is called a
$\beta$-unstable fiber} passing through $x\in H$ with $\beta\in
\mathbb{R}^+$,  if
$\|\phi(t,\omega,x)-\phi(t,\omega,\widetilde{x})\|_H=O(e^{\beta t}),
\forall\; \omega\in \Omega$ as $t\to -\infty$ for all $x,
\widetilde{x}\in \mathcal{W}_{\beta u}$.
\par
Stable fibers form a stable foliation, while unstable fibers form an unstable foliation.
Occasionally we use  $\mathcal{W}_\beta $ to denote either fibers.  Furthermore, we say a foliation is \emph{invariant } if the random
dynamical system $\phi$ maps one fiber to another fiber in the
following sense
$$
\phi(t, \omega, \mathcal{W}_\beta(x, \omega))\subset
\mathcal{W}_\beta(\phi(t,\omega, x), \theta_t\omega).
$$

\subsection{A slow-fast random dynamical system}
\label{slow888}

\quad\; Let $\Theta_1=(\Omega_1,\mathcal{F}_1, \mathds{P}_1,
\theta^1)$ and $\Theta_2=(\Omega_2,\mathcal{F}_2, \mathds{P}_2,
\theta^2)$ be two independent metric dynamical systems as introduced
in Section 2.2.
Define
$$
\Theta:=\Theta_1\times \Theta_2=(\Omega_1\times\Omega_2,
\mathcal{F}_1\otimes\mathcal{F}_2, \mathds{P}_1\times \mathds{P}_2,
(\theta^1, \theta^2)^T),
$$
and
$$
\theta_t \omega:=(\theta_t^1\omega_1, \theta_t^2\omega_2)^T,\quad
\hbox{for}\quad \omega:=(\omega_1,\omega_2)^T \in
\Omega_1\times\Omega_2:=\Omega.
$$
\par

\par
Let $W_1(t)$ and $W_2(t)$ be two mutually independent standard Wiener
processes with values in $H_s$ and $H_f$, with covariances
$Q_1=\hbox{Id}_{H_s}$ and $Q_2=\hbox{Id}_{H_f}$, respectively.

\par
Consider the following linear stochastic evolutionary equations
\begin{equation}\label{OU-3}
d\delta(t)=A\delta dt +\sigma_1 dW_1,
\end{equation}

\begin{equation}\label{OU-1}
d\eta^{\frac{1}{\varepsilon}}(t)=\frac{1}{\varepsilon}B\eta^{\frac{1}{\varepsilon}}dt
+\frac{\sigma_2}{\sqrt{\varepsilon}}dW_2,
\end{equation}
and
\begin{equation}\label{OU-2}
d\xi(t)=B\xi dt +\sigma_2 dW_2.
\end{equation}
\par
{\bf Lemma 2.1}$^{\cite{SS}}$ \quad {\it Assume that the Hypothesis
H1 holds. Then equations (\ref{OU-3}), (\ref{OU-1}) and (\ref{OU-2})
have continuous stationary solutions $\delta(\theta_t^1\omega_1)$,
$\eta^{\frac{1}{\varepsilon}}(\theta_t^{2}\omega_2)$ and
$\xi(\theta_t^2\omega_2)$, respectively. Furthermore, the stochastic
process $\eta^{\frac{1}{\varepsilon}}(\theta_t^{2}\omega_2)$ has the
same distribution as the process
$\xi(\theta_{\frac{t}{\varepsilon}}^2\omega_2)$.}
\par
Introduce new variables
\begin{equation}\label{Transformation}
X^\varepsilon=x^\varepsilon-\delta(\theta_t^1\omega_1),\quad
\hbox{and}\quad
Y^\varepsilon=y^\varepsilon-\eta^{\frac{1}{\varepsilon}}(\theta_t^2\omega_2).
\end{equation}
Then the slow-fast stochastic evolutionary equations
(\ref{Equation-s})-(\ref{Equation-f}) can be rewritten as
\begin{equation}\label{RE-s}
\frac{dX^\varepsilon}{dt}=AX^\varepsilon+F(X^\varepsilon,
Y^\varepsilon,\theta_t^\varepsilon\omega),~~~~
\end{equation}
\begin{equation}\label{RE-f}
\frac{dY^\varepsilon}{dt}=\frac{1}{\varepsilon}BY^\varepsilon+\frac{1}{\varepsilon}G(X^\varepsilon,
Y^\varepsilon,\theta_t^\varepsilon\omega),
\end{equation}
where $$F(X^\varepsilon,
Y^\varepsilon,\theta_t^\varepsilon\omega):=f(X^\varepsilon+\delta(\theta_t^1\omega_1),
Y^\varepsilon+\eta^{\frac{1}{\varepsilon}}(\theta_t^2\omega_2)),$$
and
$$G(X^\varepsilon,
Y^\varepsilon,\theta_t^\varepsilon\omega):=g(X^\varepsilon+\delta(\theta_t^1\omega_1),
Y^\varepsilon+\eta^{\frac{1}{\varepsilon}}(\theta_t^2\omega_2)).$$
The state space for this system is $H =H_s \times H_f$.
\par
Supplement the initial condition
\begin{equation}\label{initial}
X^\varepsilon(0)=X_0, \quad\hbox{and}\quad Y^\varepsilon(0)=Y_0.
\end{equation}
\par
Under the Hypothesis H1-H3, by the classical evolutionary equation
theory, system (\ref{RE-s})-(\ref{initial}) has a unique global
solution for every $\omega=(\omega_1,\omega_2)^T \in
\Omega=\Omega_1\times\Omega_2$. No exceptional sets with respect to
the initial conditions appear. Hence the solution mapping
$$ (t,\omega, (X_0,
Y_0)^T)\mapsto \Phi^\varepsilon(t,\omega, (X_0,
Y_0)^T):=(X^\varepsilon(t, \omega, (X_0, Y_0)^T), Y^\varepsilon(t,
\omega, (X_0, Y_0)^T))^T
$$
generates a continuous random dynamical system. In fact, the mapping
$\Phi^\varepsilon$ is $(\mathcal{B}(\mathbb{R})\otimes
(\mathcal{F}_1\otimes\mathcal{F}_2)\otimes \mathcal{B}(H_s\times
H_f), (\mathcal{F}_1\otimes\mathcal{F}_2))$-measurable.
\par

As in Jones \cite[p.49]{Jones-book-1995},    a \emph{slow} manifold
of a slow-fast system is the manifold in which the fast variable is
represented by the slow variable, when the scale parameter
 $\varepsilon$ is sufficiently small. It also exponentially
attracts other dynamical orbits. A \emph{critical} manifold of a
slow-fast system is the slow manifold  corresponding to the
zero scale parameter.


\par

For (\ref{RE-s})-(\ref{RE-f}), similar to  Fu, Liu and Duan
\cite{FLD} or Wang and Roberts \cite{WR-2013-JMAA}, we have the following result about  the slow manifold.
\par
Consider the so-called Liapunov-Perron equation
\begin{equation}\label{manifold-h}
h^\varepsilon(\zeta,\omega)=\frac{1}{\varepsilon}\int_{-\infty}^0e^{-\frac{Bs}{\varepsilon}}G(X^\varepsilon(s,
\omega; \zeta),  Y^\varepsilon(s, \omega; \zeta),
\theta_s^\varepsilon\omega)ds,\quad \hbox{for any}\quad \zeta\in
H_s,
\end{equation}
where $X^\varepsilon(t, \omega; \zeta)$ and $Y^\varepsilon(t,
\omega; \zeta)$ are the solutions of system
(\ref{RE-s})-(\ref{RE-f}) with the following forms
$$
\left(
\begin{array}{l}
X^\varepsilon(t, \omega; \zeta)\\
Y^\varepsilon(t, \omega; \zeta)
\end{array}
\right) =\left(
\begin{array}{l}
e^{At}\zeta+\int_0^te^{A(t-s)}F(X^\varepsilon(s, \omega;
\zeta), Y^\varepsilon(s, \omega; \zeta), \theta_s^\varepsilon\omega)ds\\
\frac{1}{\varepsilon}\int_{-\infty}^t
e^{\frac{B(t-s)}{\varepsilon}}G(X^\varepsilon(s, \omega; \zeta),
Y^\varepsilon(s, \omega; \zeta),\theta_s^\varepsilon\omega)ds
\end{array}
\right).
$$
Then
\begin{equation}\label{manifold}
\mathcal{M}^\varepsilon(\omega)=\{(\zeta,
h^\varepsilon(\zeta,\omega))^T |\; \zeta \in H_s \}
\end{equation}
is the slow manifold of the system (\ref{RE-s})-(\ref{RE-f}). It
is invariant in the following sense
$$
\Phi^\varepsilon(t,\omega, \mathcal{M}^\varepsilon(\omega))\subset
\mathcal{M}^\varepsilon(\theta_t\omega),\quad \hbox{for}\quad t > 0.
$$
Furthermore, the slow manifold exponentially attracts other
dynamical orbits.


In the rest of this paper, we   use $\mathcal{W}_\beta^\varepsilon((X_0,
Y_0)^T, \omega)$ to denote a fiber of the slow foliation, and use
$\mathcal{W}_\beta^0((X_0, Y_0)^T, \omega)$ to denote a fiber of the
critical foliation. According to Section 2.2, the slow foliation is
essentially   an unstable foliation.
\par

\renewcommand{\theequation}{\thesection.\arabic{equation}}
\setcounter{equation}{0}

\section{A motivating example for slow foliation}

\quad\; Before   presenting a general theory, we work out a simple   example   for   slow foliation.
\par

Consider the following slow-fast stochastic ordinary differential
equations
\begin{equation}\label{Equation-s-motivation}
\frac{dx^\varepsilon}{dt}=x^\varepsilon,~~~~~~~~~~~~~~~~~~~~~~~~~~~~~\quad
\hbox{in} \quad \mathbb{R}^1,
\end{equation}
\begin{equation}\label{Equation-f-motivation}
\frac{dy^\varepsilon}{dt}=-\frac{1}{\varepsilon}y^\varepsilon
+\frac{1}{\varepsilon}(x^\varepsilon)^2
+\frac{1}{\sqrt{\varepsilon}}\dot{W_2},\quad \hbox{in} \quad
\mathbb{R}^1,
\end{equation}
where   $W_2$ is a scalar Wiener process.
It follows from \S \ref{slow888} that the converted random system is as follows
\par
\begin{equation}\label{Random-s-motivation}
\frac{dX^\varepsilon}{dt}=X^\varepsilon,~~~~~~~~~~~~~~~~~~~\quad
\hbox{in} \quad \mathbb{R}^1,
\end{equation}
\begin{equation}\label{Random-f-motivation}
\frac{dY^\varepsilon}{dt}=-\frac{1}{\varepsilon}Y^\varepsilon
+\frac{1}{\varepsilon}(X^\varepsilon)^2 ,\quad \hbox{in} \quad
\mathbb{R}^1.
\end{equation}
With the initial condition $X^\varepsilon(0)=X_0$ and $Y^\varepsilon(0)=Y_0$, the solution is
\begin{equation}
X^\varepsilon(t)=X_0 e^t,~~~~~~~~~~~~~~~~~~~~~~~~~~~~~~~~~~~~t\in
\mathbb{R},
\end{equation}
\begin{equation}
Y^\varepsilon(t)=Y_0 e^{-\frac{t}{\varepsilon}}+\frac{1}{1+2\varepsilon}
X_0^2[e^{2t}-e^{-\frac{t}{\varepsilon}}],\quad t\in \mathbb{R},
\end{equation}
where
$$
\begin{array}{l}
X^\varepsilon(t)=X^\varepsilon(t,\omega_2, (X_0, Y_0)^T)=X^\varepsilon(t,\eta^{\frac{1}{\varepsilon}}(\theta_t^2\omega_2), (X_0, Y_0)^T)=x^\varepsilon,\\
Y^\varepsilon(t)=Y^\varepsilon(t,\omega_2, (X_0,
Y_0)^T)=Y^\varepsilon(t,\eta^{\frac{1}{\varepsilon}}(\theta_t^2\omega_2),
(X_0,
Y_0)^T)=y^\varepsilon-\eta^{\frac{1}{\varepsilon}}(\theta_t^2\omega_2).
\end{array}
$$
\par
For  every two points $(X_0, Y_0)^T$ and $(\widetilde{X}_0,
\widetilde{Y}_0)^T$ in $\mathbb{R}^1\times\mathbb{R}^1$, we calculate the difference between two orbits
$$
\begin{array}{ll}
I&:=|(X^\varepsilon(t,\omega_2, (X_0, Y_0)^T),
Y^\varepsilon(t,\omega_2, (X_0,
Y_0)^T))^T-(X^\varepsilon(t,\omega_2, (\widetilde{X}_0,
\widetilde{Y}_0)^T), Y^\varepsilon(t,\omega_2,
(\widetilde{X}_0, \widetilde{Y}_0)^T))^T|\\
&=|X^\varepsilon(t,\omega_2, (X_0, Y_0)^T)-X^\varepsilon(t,\omega_2,
(\widetilde{X}_0, \widetilde{Y}_0)^T)|+|Y^\varepsilon(t,\omega_2,
(X_0, Y_0)^T)-Y^\varepsilon(t,\omega_2,
(\widetilde{X}_0, \widetilde{Y}_0)^T)|\\
&\leq |X_0-\widetilde{X}_0|\cdot e^t
+\frac{1}{1+2\varepsilon}|(X_0^2-\widetilde{X}_0^2)|\cdot e^{2t}
+|(Y_0-\widetilde{Y}_0) -\frac{1}{1+2\varepsilon}(X_0^2
-\widetilde{X}_0^2)|\cdot e^{-\frac{t}{\varepsilon}}.
\end{array}
$$
If the coefficient
\begin{equation}\label{condition}
(Y_0-\widetilde{Y}_0) -\frac{1}{1+2\varepsilon}(X_0^2
-\widetilde{X}_0^2)=0,
\end{equation}
then the difference of two   orbits is
$$
I=O(e^{t}), \quad \hbox{as}\quad t\to -\infty.
$$
\par
Define
\begin{equation}\label{foliation-motivation}
\mathcal{W}_1^\varepsilon((X_0, Y_0)^T, \omega_2)=\{(\zeta,
l^\varepsilon(\zeta, (X_0, Y_0)^T, \omega_2))^T |\; \zeta\in
\mathbb{R}^1\},
\end{equation}
where the function
\begin{equation}\label{foliation-l-motivation}
l^\varepsilon(\zeta, (X_0, Y_0)^T,
\omega_2)=Y_0+\frac{1}{1+2\varepsilon}(\zeta^2-X_0^2),\quad \zeta\in
\mathbb{R}^1.
\end{equation}
\par
Whenever an initial point  $(\widetilde{X}_0, \widetilde{Y}_0)^T$ is
in $\mathcal{W}_1^\varepsilon((X_0, Y_0)^T, \omega_2)$, the
condition (\ref{condition}) holds. This immediately implies that the
different dynamical orbits will be exponentially approaching each
other as $t\to -\infty$. Therefore, we   say that
$\mathcal{W}_1^\varepsilon((X_0, Y_0)^T, \omega_2)$ is \emph{ a
fiber of the slow foliation}. It is the graph of
$l^\varepsilon(\zeta, (X_0, Y_0)^T, \omega_2)$.   Different
  orbits of the slow-fast system
(\ref{Random-s-motivation})-(\ref{Random-f-motivation}) are
exponentially approaching each other in backward time, whenever they start from the same
fiber.
\par
As seen in (\ref{foliation-l-motivation}), the slow foliation of
(\ref{Random-s-motivation})-(\ref{Random-f-motivation}) is a family
of the parallel parabolic curves (i.e., fibers) in the state space.
\par
In addition, from (\ref{manifold-h}) and (\ref{manifold}), we know
that the \emph{slow manifold} of
(\ref{Random-s-motivation})-(\ref{Random-f-motivation}) is
\begin{equation}\label{manifold-motivation}
\mathcal{M}^\varepsilon(\omega_2)=\{(\zeta,
 h^\varepsilon(\zeta,\omega_2))^T |\; \zeta\in \mathbb{R}^1\},
\end{equation}
where
\begin{equation}\label{manifold-h-motivation}
h^\varepsilon(\zeta,\omega_2)=\frac{1}{1+2\varepsilon}\zeta^2,\quad
\zeta\in \mathbb{R}^1.
\end{equation}
By comparing with (\ref{foliation-l-motivation}), it is clear that the slow manifold is a fiber of the slow
foliation.
\par

Now we consider another stochastic system independent of
$\varepsilon$ as follows
\begin{equation}\label{Equation-s-critical-example}
\frac{dx^0}{dt}=0,~~~~~~~~~~~~~~~~~~~~~~\quad \hbox{in} \quad
\mathbb{R}^1,
\end{equation}
\begin{equation}\label{Equation-f-critical-example}
\frac{dy^0}{dt}=-y^0 +(x^0)^2+\dot{W_2},\quad \hbox{in} \quad
\mathbb{R}^1.
\end{equation}
It follows from \S \ref{slow888} that the converted  random system is
\begin{equation}\label{Equation-s-R-critical-example}
\frac{dX^0}{dt}=0,~~~~~~~~~~~~~~~~~~\quad \hbox{in} \quad
\mathbb{R}^1,
\end{equation}
\begin{equation}\label{Equation-f-R-critical-example}
\frac{dY^0}{dt}=-Y^0 +(X^0)^2,\quad \hbox{in} \quad \mathbb{R}^1.
\end{equation}
The solution  with   initial condition $X^0(0)=X_0$ and $Y^0(0)=Y_0$ is
$$
\begin{array}{ll}
X^0(t)=X_0,\quad & t\in \mathbb{R},\\
Y^0(t)=e^{-t}Y_0+X_0^2(1-e^{-t}),\quad & t\in \mathbb{R},
\end{array}
$$
where
$$
\begin{array}{l}
X^0(t)=X^0(t,\omega_2, (X_0, Y_0)^T)=X^0(t,\xi(\theta_t^2\omega_2),
(X_0,
Y_0)^T)=x^0,\\
Y^0(t)=Y^0(t,\omega_2, (X_0, Y_0)^T)=Y^0(t,\xi(\theta_t^2\omega_2),
(X_0, Y_0)^T)=y^0-\xi(\theta_t^2\omega_2).
\end{array}
$$
\par
By the same argument as   above,   a fiber of the foliation of
(\ref{Equation-s-R-critical-example})-(\ref{Equation-f-R-critical-example})
is
\begin{equation}
\mathcal{W}_1^0((X_0, Y_0)^T, \omega_2)=\{(\zeta, l^0(\zeta, (X_0,
Y_0)^T, \omega_2))^T |\; \zeta\in \mathbb{R}^1\},
\end{equation}
where
\begin{equation}\label{Foliation-l-example}
l^0(\zeta, (X_0, Y_0)^T, \omega_2)=Y_0+(\zeta^2-X_0^2),\quad
\zeta\in \mathbb{R}^1.
\end{equation}}
This  is called the \emph{critical foliation} for the system
(\ref{Random-s-motivation})-(\ref{Random-f-motivation}).
\par
Observe that, by a time change $t\to \varepsilon t$,   Equation
(\ref{Random-f-motivation}) is transformed to Equation
(\ref{Equation-f-R-critical-example}). Also notice that
$\eta^{\frac{1}{\varepsilon}}(\theta_{\varepsilon t}^2\omega_2)$ has
the same distribution as  $\xi(\theta_{t}^2\omega_2)$ by Lemma 2.1.
Thus
$$
l^\varepsilon(\zeta, (X_0, Y_0)^T,
\omega_2)=Y_0+\frac{1}{1+2\varepsilon}(\zeta^2-X_0^2)
\overset{\text{d}}{\longrightarrow } Y_0+(\zeta^2-X_0^2)=l^0(\zeta,
(X_0, Y_0)^T, \omega_2),\quad \hbox{as}\quad \varepsilon\to 0,
$$
for $\zeta\in \mathbb{R}^1$, where
``$\overset{\text{d}}{\longrightarrow }$" denotes the convergence in
distribution.
Therefore, the slow foliation converges in distribution to the critical foliation, as $\varepsilon$ tends to zero.
\par

\renewcommand{\theequation}{\thesection.\arabic{equation}}
\setcounter{equation}{0}

\section{Slow foliation}

\quad\quad In this section, we establish a   theory of the slow
foliation for the slow-fast system
(\ref{Equation-s})-(\ref{Equation-f}). We derive the existence of
slow foliation for the corresponding  random slow-fast system
(\ref{RE-s})-(\ref{RE-f}). The dynamical orbits of the system
(\ref{RE-s})-(\ref{RE-f}) in a given fiber are shown to
exponentially approach each other   in backward time. In addition,
we explore the geometric structure of the slow foliation and analyze
the relationship between the slow foliation and the slow manifold.
\par

Define two Banach spaces for a fixed $\beta$ as follows:
$$
C_{\beta}^{s,-}=\{\varphi:(-\infty, 0]\to H_s|\quad
\varphi\;\hbox{is
 continous and}\; \sup\limits_{t\leq 0}e^{-\beta t}
\|\varphi\|_{H_s}<\infty \},
$$
$$
C_{\beta}^{f,-}=\{\varphi:(-\infty, 0]\to H_f|\quad
\varphi\;\hbox{is
 continous and}\; \sup\limits_{t\leq 0}e^{-\beta t}
\|\varphi\|_{H_f}<\infty \},
$$
with the norms
$$\|\varphi\|_{C_{\beta}^{s,-}}=\sup\limits_{t\leq
0}e^{-\beta t} \|\varphi\|_{H_s},\quad
\hbox{and}\quad\|\varphi\|_{C_{\beta}^{f,-}}=\sup\limits_{t\leq
0}e^{-\beta t} \|\varphi\|_{H_f},
$$
respectively. Define $C_\beta^-:=C_\beta^{s,-}\times C_\beta^{f,-}$,
with norm   $\|(X,
Y)^T\|_{C_\beta^-}=\|X\|_{C_\beta^{s,-}}+\|Y\|_{C_\beta^{f,-}}$,  for
$(X, Y)^T\in C_\beta^-$.
\par

Denote $\Phi^\varepsilon(t,\omega,
(X_0,Y_0)^T)=(X^\varepsilon(t,\omega, (X_0,Y_0)^T),
Y^\varepsilon(t,\omega, (X_0,Y_0)^T))^T$ the solution of the
slow-fast random system (\ref{RE-s})-(\ref{RE-f}) with the initial
condition $\Phi^\varepsilon(0,\omega, (X_0,Y_0)^T)=(X_0,Y_0)^T$.
\par
Define the difference of two dynamical orbits
\begin{equation}\label{Equation-Psi}
\begin{array}{ll}
\Psi^\varepsilon(t)&=\Phi^\varepsilon(t,\omega,
(\widetilde{X}_0,\widetilde{Y}_0)^T)-\Phi^\varepsilon(t,\omega,
(X_0,Y_0)^T)\\
&=(X^\varepsilon(t,\omega,
(\widetilde{X}_0,\widetilde{Y}_0)^T)-X^\varepsilon(t,\omega,
(X_0,Y_0)^T), Y^\varepsilon(t,\omega,
(\widetilde{X}_0,\widetilde{Y}_0)^T)-Y^\varepsilon(t,\omega,
(X_0,Y_0)^T))^T\\
&:=(U^\varepsilon(t), V^\varepsilon(t))^T.
\end{array}
\end{equation}
Here the initial condition
$$
\Psi^\varepsilon(0)=(U^\varepsilon(0),
V^\varepsilon(0))^T=(\widetilde{X}_0-X_0, \widetilde{Y}_0-Y_0)^T,
$$
and the solution
$$
\begin{array}{l}
X^\varepsilon(t,\omega,
(\widetilde{X}_0,\widetilde{Y}_0)^T)=U^\varepsilon(t)+X^\varepsilon(t,\omega,
(X_0,Y_0)^T),\\
Y^\varepsilon(t,\omega,
(\widetilde{X}_0,\widetilde{Y}_0)^T)=V^\varepsilon(t)+Y^\varepsilon(t,\omega,
(X_0,Y_0)^T).
\end{array}
$$
Moreover, $(U^\varepsilon, V^\varepsilon)^T$ satisfies
\begin{equation}\label{Equation-U}
\frac{dU^\varepsilon}{dt}=AU^\varepsilon+\Delta F(U^\varepsilon,
V^\varepsilon, \theta_t^\varepsilon\omega),~~~~
\end{equation}
\begin{equation}\label{Equation-V}
\frac{dV^\varepsilon}{dt}=\frac{1}{\varepsilon}BV^\varepsilon+\frac{1}{\varepsilon}\Delta
G(U^\varepsilon, V^\varepsilon, \theta_t^\varepsilon\omega),
\end{equation}
with    nonlinearities
\begin{equation}\label{Equation-DF}
\begin{array}{lll}
\Delta F(U^\varepsilon, V^\varepsilon, \theta_t^\varepsilon\omega)
&=& F(U^\varepsilon(t)+X^\varepsilon(t,\omega, (X_0,Y_0)^T),
V^\varepsilon(t)+Y^\varepsilon(t,\omega, (X_0,Y_0)^T),
\theta_t^\varepsilon\omega)\\
&&-F(X^\varepsilon(t,\omega, (X_0,Y_0)^T), Y^\varepsilon(t,\omega,
(X_0,Y_0)^T), \theta_t^\varepsilon\omega),
\end{array}
\end{equation}
and
\begin{equation}\label{Equation-DG}
\begin{array}{lll}
\Delta G(U^\varepsilon, V^\varepsilon, \theta_t^\varepsilon\omega)
&=& G(U^\varepsilon(t)+X^\varepsilon(t,\omega, (X_0,Y_0)^T),
V^\varepsilon(t)+Y^\varepsilon(t,\omega, (X_0,Y_0)^T),
\theta_t^\varepsilon\omega)\\
&&-G(X^\varepsilon(t,\omega, (X_0,Y_0)^T), Y^\varepsilon(t,\omega,
(X_0,Y_0)^T), \theta_t^\varepsilon\omega),
\end{array}
\end{equation}
and   initial condition
$$
U^\varepsilon(0)=U_0=\widetilde{X}_0-X_0,
V^\varepsilon(0)=V_0=\widetilde{Y}_0-Y_0.
$$
\par

Define
\begin{equation}\label{Foliation}
\mathcal{W}_\beta^\varepsilon((X_0, Y_0)^T,
\omega)=\{(\widetilde{X}_0, \widetilde{Y}_0)^T\in H_s\times H_f| \;
\Phi^\varepsilon(t,\omega, (X_0,Y_0)^T)-\Phi^\varepsilon(t,\omega,
(\widetilde{X}_0,\widetilde{Y}_0)^T)\in C_{\beta}^-\}.
\end{equation}
As we will show, $\mathcal{W}_\beta^\varepsilon((X_0, Y_0)^T,
\omega)$ is a fiber of the slow foliation for the slow-fast random
system (\ref{RE-s})-(\ref{RE-f}).
\par
Now we present some lemmas before   our main results.
\par

{\bf Lemma 4.1}\quad {\it Assume that the Hypotheses H1-H3 hold.
Take $\beta$ as the positive real number $\frac{-\gamma_s}{2}$. Then
$(\widetilde{X}_0, \widetilde{Y}_0)^T$ is in
$\mathcal{W}_\beta^\varepsilon((X_0, Y_0)^T, \omega)$ if and only if
there exists a function
$\Psi^\varepsilon(t)=(U^\varepsilon(t),
V^\varepsilon(t))^T=(U^\varepsilon(t, \omega, (X_0, Y_0)^T ;
U^\varepsilon(0)), V^\varepsilon(t, \omega, (X_0, Y_0)^T ;
U^\varepsilon(0)))^T\in C_\beta^-$   such that
\begin{equation}\label{Lemma 3.1-0}
\Psi^\varepsilon(t)=\left(
\begin{array}{l}
U^\varepsilon(t)\\
V^\varepsilon(t)
\end{array}
\right) = \left(
\begin{array}{c}
e^{At}U^\varepsilon(0)+\int_0^te^{A(t-s)}\Delta F(U^\varepsilon(s),
V^\varepsilon(s), \theta_s^\varepsilon\omega)ds
\\
\frac{1}{\varepsilon}\int_{-\infty}^te^{\frac{B(t-s)}{\varepsilon}}\Delta
G(U^\varepsilon(s), V^\varepsilon(s), \theta_s^\varepsilon\omega)ds
\end{array}
\right),
\end{equation}
where $\Delta F$ and $\Delta G$ are defined in (\ref{Equation-DF})
and (\ref{Equation-DG}). }
\par

{\bf Proof.}\quad   Let $(\widetilde{X}_0, \widetilde{Y}_0)^T\in
\mathcal{W}_\beta^\varepsilon((X_0, Y_0)^T,\omega)$. By the
variation of constants formula,  we have
\begin{equation}\label{Lemma 3.1-1}
U^\varepsilon(t)=e^{At}U^\varepsilon(0)+\int_0^te^{A(t-s)}\Delta
F(U^\varepsilon(s), V^\varepsilon(s),
\theta_s^\varepsilon\omega)ds,~~~~~~~
\end{equation}
\begin{equation}\label{Lemma 3.1-2}
V^\varepsilon(t)=e^{\frac{B(t-\tau)}{\varepsilon}}V^\varepsilon(\tau)
+\frac{1}{\varepsilon}\int_{\tau}^te^{\frac{B(t-s)}{\varepsilon}}\Delta
G(U^\varepsilon(s), V^\varepsilon(s), \theta_s^\varepsilon\omega)ds.
\end{equation}
Since $\Phi^\varepsilon(\cdot)\in C_\beta^-$, for $\tau<0$,
$$
\begin{array}{ll}
\|e^{\frac{B(t-\tau)}{\varepsilon}}V^\varepsilon(\tau)\|_{H_f}&\leq
\|V\|_{C_\beta^{f,-}}\cdot e^{\frac{-\gamma_f\cdot
t}{\varepsilon}}\cdot
e^{(\frac{-\gamma_f}{\varepsilon}-\beta)\cdot(-\tau)}\\
&\leq \|V\|_{C_\beta^{f,-}}\cdot e^{\frac{-\gamma_f\cdot
t}{\varepsilon}}\cdot
e^{(\frac{-\gamma_f}{\varepsilon}-\beta)\cdot(-\tau)}\\
& \longrightarrow 0,\quad\hbox{as}\quad \tau\to -\infty,
\end{array}
$$
which implies that
\begin{equation}\label{Lemma 3.1-3}
V^\varepsilon(t)=\frac{1}{\varepsilon}\int_{-\infty}^te^{\frac{B(t-s)}{\varepsilon}}\Delta
G(U^\varepsilon(s), V^\varepsilon(s), \theta_s^\varepsilon\omega)ds.
\end{equation}
\par
Therefore, it follows from (\ref{Lemma 3.1-1})-(\ref{Lemma 3.1-3})
that (\ref{Lemma 3.1-0}) holds.  By direct calculation,
it   is clear that the converse holds. This completes the proof of Lemma 4.1.
\hfill{$\blacksquare$}
\par

{\bf Lemma 4.2}\quad {\it Assume that the Hypotheses  H1-H3 hold.
Take $\beta$ as the positive real number $\frac{-\gamma_s}{2}$. For
any given $U^\varepsilon (0)=\widetilde{X}_0-X_0\in H_s$, there
exists a sufficiently small positive constant $\varepsilon_0$ such
that for $\varepsilon\in (0, \varepsilon_0)$, the system (\ref{Lemma
3.1-0}) has a unique solution
$\Psi^\varepsilon(\cdot)=\Psi^\varepsilon(\cdot, \omega, (X_0,
Y_0)^T ; U^\varepsilon(0))$ in $C_{\beta}^-$.}
\par
{\bf Proof.}\quad Introduce two operators
$\mathcal{J}_s^\varepsilon: C_\beta^-\longrightarrow C_\beta^{s,-}$
and $\mathcal{J}_f^\varepsilon: C_\beta^-\longrightarrow
C_\beta^{f,-}$ satisfying
$$
\mathcal{J}_s^\varepsilon(\Psi^\varepsilon)[t]
=e^{At}U^\varepsilon(0)+\int_0^te^{A(t-s)}\Delta F(U^\varepsilon(s),
V^\varepsilon(s), \theta_s^\varepsilon\omega)ds,
$$
$$
\mathcal{J}_f^\varepsilon(\Psi^\varepsilon)[t]
=\frac{1}{\varepsilon}\int_{-\infty}^te^{\frac{B(t-s)}{\varepsilon}}\Delta
G(U^\varepsilon(s), V^\varepsilon(s),
\theta_s^\varepsilon\omega)ds.~~~~~~~~~~
$$
Also pose the operator $\mathcal{J}^\varepsilon:
C_\beta^-\longrightarrow C_\beta^{-}$ satisfying
$\mathcal{J}^\varepsilon(\Psi^\varepsilon)=
(\mathcal{J}_s^\varepsilon(\Psi^\varepsilon),
\mathcal{J}_f^\varepsilon(\Psi^\varepsilon))^T$. It is easy to
verify that $\mathcal{J}_s^\varepsilon$, $\mathcal{J}_f^\varepsilon$
and $\mathcal{J}^\varepsilon$ are well-defined in $C_{\beta}^{s,-}$,
$C_{\beta}^{f,-}$ and $C_{\beta}^-$, respectively.
\par

For any $\Psi^\varepsilon=(U^\varepsilon, V^\varepsilon)^T\in
C_\beta^-$ and
$\widetilde{\Psi}^\varepsilon=(\widetilde{U}^\varepsilon,
\widetilde{V}^\varepsilon)^T\in C_\beta^-$, then
\begin{equation}\label{Lemma 3.2-1}
\begin{array}{ll}
&\|\mathcal{J}_s^\varepsilon(\Psi^\varepsilon)- \mathcal{J}_s^\varepsilon(\widetilde{\Psi}^\varepsilon)\|_{C_\beta^{s,-}}\\
=&\|\int_{0}^t e^{A(t-s)}[\Delta F(U^\varepsilon(s),
V^\varepsilon(s),\theta_s^\varepsilon\omega)- \Delta
F(\widetilde{U}^\varepsilon(s),
\widetilde{V}^\varepsilon(s),\theta_s^\varepsilon\omega)]ds
\|_{C_\beta^{s,-}}\\
=&\|\int_{0}^te^{A(t-s)}[F(U^\varepsilon(s)+X^\varepsilon(s,\omega,
(X_0, Y_0)^T), V^\varepsilon(s)+Y^\varepsilon(s,\omega,
(X_0, Y_0)^T),\theta_s^\varepsilon\omega)\\
&~~~~~- F(\widetilde{U}^\varepsilon(s)+X^\varepsilon(s,\omega, (X_0,
Y_0)^T), \widetilde{V}^\varepsilon(s)+Y^\varepsilon(s,\omega, (X_0,
Y_0)^T),\theta_s^\varepsilon\omega)]ds
\|_{C_\beta^{s,-}}\\
\leq &\sup\limits_{t\leq 0}e^{-\beta t}\cdot
K\int_{t}^0e^{-\gamma_s(t-s)}
(\|U^\varepsilon(s)-\widetilde{U}^\varepsilon(s)\|_{H_s}+\|V^\varepsilon(s)-\widetilde{V}^\varepsilon(s)\|_{H_f})ds\\
\leq
&\frac{K}{-\beta-\gamma_s}\|\Psi^\varepsilon-\widetilde{\Psi}^\varepsilon\|_{C_\beta^-},
\end{array}
\end{equation}
and
\begin{equation}\label{Lemma 3.2-2}
\begin{array}{ll}
&\|\mathcal{J}_f^\varepsilon(\Psi^\varepsilon)- \mathcal{J}_f^\varepsilon(\widetilde{\Psi}^\varepsilon)\|_{C_\beta^{f,-}}\\
=&\|\frac{1}{\varepsilon}\int_{-\infty}^te^{\frac{B(t-s)}{\varepsilon}}[\Delta
G(U^\varepsilon(s), V^\varepsilon(s),\theta_s^\varepsilon\omega)-
\Delta G(\widetilde{U}^\varepsilon(s),
\widetilde{V}^\varepsilon(s),\theta_s^\varepsilon\omega)]ds
\|_{C_\beta^{f,-}}\\
=&\|\frac{1}{\varepsilon}\int_{-\infty}^te^{\frac{B(t-s)}{\varepsilon}}[G(U^\varepsilon(s)+X^\varepsilon(s,\omega,
(X_0, Y_0)^T), V^\varepsilon(s)+Y^\varepsilon(s,\omega,
(X_0, Y_0)^T),\theta_s^\varepsilon\omega)\\
&~~~~~- G(\widetilde{U}^\varepsilon(s)+X^\varepsilon(s,\omega, (X_0,
Y_0)^T), \widetilde{V}^\varepsilon(s)+Y^\varepsilon(s,\omega, (X_0,
Y_0)^T),\theta_s^\varepsilon\omega)]ds
\|_{C_\beta^{f,-}}\\
\leq &\sup\limits_{t\leq 0}e^{-\beta
t}\cdot\frac{K}{\varepsilon}\int_{-\infty}^te^{\frac{-\gamma_f(t-s)}{\varepsilon}}
(\|U^\varepsilon(s)-\widetilde{U}^\varepsilon(s)\|_{H_s}+\|V^\varepsilon(s)-\widetilde{V}^\varepsilon(s)\|_{H_f})ds\\
\leq &\frac{K}{\gamma_f+\varepsilon
\beta}\|\Psi^\varepsilon-\widetilde{\Psi}^\varepsilon\|_{C_\beta^-}.
\end{array}
\end{equation}
\par
It immediately follows from (\ref{Lemma 3.2-1}) and (\ref{Lemma
3.2-2}) that
\begin{equation}\label{Lemma 3.2-3}
\begin{array}{ll}
\|\mathcal{J}^\varepsilon(\Psi^\varepsilon)-
\mathcal{J}^\varepsilon(\widetilde{\Psi}^\varepsilon)\|_{C_\beta^-}
&=\|\mathcal{J}_s^\varepsilon(\Psi^\varepsilon)-
\mathcal{J}_s^\varepsilon(\widetilde{\Psi}^\varepsilon)\|_{C_\beta^{s,-}}+\|\mathcal{J}_f^\varepsilon(\Psi^\varepsilon)-
\mathcal{J}_f^\varepsilon(\widetilde{\Psi}^\varepsilon)\|_{C_\beta^{f,-}}\\
&\leq (\frac{K}{-\beta-\gamma_s}+\frac{K}{\gamma_f+\varepsilon
\beta})\|\Psi^\varepsilon-\widetilde{\Psi}^\varepsilon\|_{C_\beta^-}.
\end{array}
\end{equation}
Put the constant
\begin{equation}\label{Lemma 3.2-4}
\rho(\gamma_s,\gamma_f,K,\varepsilon)=\frac{K}{-\beta-\gamma_s}+\frac{K}{\gamma_f+\varepsilon
\beta}.
\end{equation}
Then
\begin{equation}\label{Lemma 3.2-5}
\|\mathcal{J}^\varepsilon(\Psi^\varepsilon)-
\mathcal{J}^\varepsilon(\widetilde{\Psi}^\varepsilon)\|_{C_\beta^-}\leq
\rho(\gamma_s,\gamma_f,K,\varepsilon)
\|\Psi^\varepsilon-\widetilde{\Psi}^\varepsilon\|_{C_\beta^-}.
\end{equation}
\par
Notice that the Hypothesis H3 holds, $\beta=\frac{-\gamma_s}{2}$,
and that
$$
\rho(\gamma_s,\gamma_f,K,\varepsilon)\longrightarrow
\frac{K}{-\beta-\gamma_s}+\frac{K}{\gamma_f},\quad \hbox{as}\quad
\varepsilon\searrow 0.
$$
Therefore, there is a sufficiently small
positive constant $\varepsilon_0$ such that for $\varepsilon\in (0,
\varepsilon_0)$,
$$
0<\rho(\gamma_s,\gamma_f,K,\varepsilon)<1.
$$
Then the map $\mathcal{J}^\varepsilon(\Psi^\varepsilon)$ is
contractive in $C_\beta^-$ uniformly with respect to $(\omega, (X_0,
Y_0)^T , U^\varepsilon(0))$. By the contraction mapping principle,
we have that for each $U^\varepsilon(0)\in H_s$, the mapping
$\mathcal{J}^\varepsilon(\Psi^\varepsilon)=\mathcal{J}^\varepsilon(\Psi^\varepsilon,
\omega, (X_0, Y_0)^T ; U^\varepsilon(0))$ has a unique fixed point,
which still denoted by $$\Psi^\varepsilon(\cdot)
=\Psi^\varepsilon(\cdot, \omega, (X_0, Y_0)^T; U^\varepsilon(0)) \in
C_\beta^-.$$ In other words, $\Psi^\varepsilon(\cdot, \omega, (X_0,
Y_0)^T ; U^\varepsilon(0))\in C_\beta^-$ is a unique solution of the
system (\ref{Lemma 3.1-0}).\hfill{$\blacksquare$}
\par

{\bf Lemma 4.3}\quad {\it Assume that the Hypothesis H1-H3 hold.
Take $\beta$ as the positive real number $\frac{-\gamma_s}{2}$. Let
$\Psi^\varepsilon(t)=\Psi^\varepsilon(t, \omega, (X_0, Y_0)^T ;
U^\varepsilon(0))$ be the unique solution of the system (\ref{Lemma
3.1-0}) in $C_{\beta}^-$. For any $U^\varepsilon(0)$ and
$\widetilde{U}^\varepsilon(0)$ in $C_{\beta}^{s,-}$, then there
exists a sufficiently small positive constant $\varepsilon_0$ such
that for $\varepsilon\in (0, \varepsilon_0)$, we have
\begin{equation}
\begin{array}{ll}
&\|\Psi^\varepsilon(t, \omega, (X_0, Y_0)^T ;
U^\varepsilon(0))-\Psi^\varepsilon(t, \omega, (X_0, Y_0)^T ; \widetilde{U}^\varepsilon(0))\|_{C_{\beta}^-}\\
\leq &
\frac{1}{1-\rho(\gamma_s,\gamma_f,K,\varepsilon)}\|U^\varepsilon(0)-\widetilde{U}^\varepsilon(0)\|_{H_s},
\end{array}
\end{equation}
where $\rho(\gamma_s,\gamma_f,K,\varepsilon)$ is defined as
(\ref{Lemma 3.2-4}). }
\par
Lemma 4.3 is easily deduced by using the same arguments as in the
proof of Lemma 4.2. Here we omit it.
\par
In the following, for every $\zeta\in H_s$, we define
\begin{equation}\label{Foliation-l}
\begin{array}{ll}
& l^\varepsilon(\zeta,(X_0, Y_0)^T,\omega)\\
:=&Y_0
+\frac{1}{\varepsilon}\int_{-\infty}^0e^{\frac{-Bs}{\varepsilon}}
\Delta G(U^\varepsilon(s, \omega, (X_0, Y_0)^T ; (\zeta-X_0)),\\
&~~~~~~~~~~~~~~~~~~~~~~~~~~~~V^\varepsilon(s, \omega, (X_0, Y_0)^T ;
(\zeta-X_0)), \theta_s^\varepsilon\omega)ds.
\end{array}
\end{equation}
\par
Now we give our main result.
\par
{\bf Theorem 4.1 (Slow foliation)}\quad{\it Assume that the
Hypothesis H1-H3 hold. Take $\beta$ as the positive real number
$\frac{-\gamma_s}{2}$. Then there exists a sufficiently small
positive constant $\varepsilon_0$ such that for $\varepsilon\in (0,
\varepsilon_0)$, the invariant foliation of the slow-fast random
system (\ref{RE-s})-(\ref{RE-f}) exists.
\par
(i)\quad Its one fiber is the graph of a Lipschitz function. That is
\begin{equation}\label{Foliation-mapping}
\mathcal{W}_\beta^\varepsilon((X_0, Y_0)^T, \omega)=\{(\zeta,
l^\varepsilon(\zeta, (X_0, Y_0)^T, \omega))^T| \; \zeta\in H_s\},
\end{equation}
where $(X_0, Y_0)^T\in H_s\times H_f$, and the function
$l^\varepsilon(\zeta,(X_0, Y_0)^T, \omega)$ is defined as
(\ref{Foliation-l}). In addition, $l^\varepsilon(\zeta,(X_0, Y_0)^T,
\omega)$ is Lipschitz continuous with respect to $\zeta$, whose
Lipschitz constant $Lip l^\varepsilon$ satisfies
$$
Lip l^\varepsilon \leq
\frac{K}{\gamma_f+\varepsilon\beta}\cdot\frac{1}{1-\rho(\gamma_s,\gamma_f,K,\varepsilon)},
$$
where $\rho(\gamma_s,\gamma_f,K,\varepsilon)$ is defined as
(\ref{Lemma 3.2-4}).
\par
(ii)\quad 
The dynamical orbits of (\ref{RE-s})-(\ref{RE-f}) are exponentially
approaching each other in backward time only if they start from the
same fiber. That is, for any two points $(\widetilde{X}_0^1,
\widetilde{Y}_0^1)^T$ and $(\widetilde{X}_0^2, \widetilde{Y}_0^2)^T$
in a same fiber $\mathcal{W}_\beta^\varepsilon((X_0, Y_0)^T,
\omega)$,
\begin{equation}\label{Foliation-attracting}
\begin{array}{ll}
\|\Phi^\varepsilon(t,\omega, (\widetilde{X}_0^1,
\widetilde{Y}_0^1)^T)-\Phi^\varepsilon(t,\omega, (\widetilde{X}_0^2,
\widetilde{Y}_0^2)^T)\|_{H_s\times H_f}&\leq \frac{e^{\beta
t}}{1-\rho(\gamma_s,\gamma_f,K,\varepsilon)}\cdot
\|\widetilde{X}_0^1-\widetilde{X}_0^2\|_{H_s}\\
& =O(e^{\beta t}),~~~~~~~~~~ \forall\quad t\to-\infty.
\end{array}
\end{equation}
\par
(iii)\quad Its fiber is invariant, i.e.,
$$
\Phi^\varepsilon(t, \omega, \mathcal{W}_\beta^\varepsilon((X_0,
Y_0)^T, \omega))\subset
\mathcal{W}_\beta^\varepsilon(\Phi^\varepsilon(t,\omega, (X_0,
Y_0)^T), \theta_t\omega).
$$
\par
}
\par
{\bf Proof.}\quad (i)\emph{ To prove a fiber of the slow foliation
is the graph of a Lipschitz function.}
\par
It follows from (\ref{Lemma 3.1-0}) that
$$
\left(
\begin{array}{l}
\widetilde{X}_0-X_0\\
\widetilde{Y}_0-Y_0
\end{array}
\right ) = \left(
\begin{array}{c}
\widetilde{X}_0-X_0\\
\frac{1}{\varepsilon}\int_{-\infty}^0e^{\frac{-Bs}{\varepsilon}}\Delta
G(U^\varepsilon(s), V^\varepsilon(s), \theta_s^\varepsilon\omega)ds
\end{array}
\right ),
$$
which implies that
$$
\begin{array}{ll}
\widetilde{Y}_0&=Y_0
+\frac{1}{\varepsilon}\int_{-\infty}^0e^{\frac{-Bs}{\varepsilon}}
\Delta G(U^\varepsilon(s, \omega, (X_0, Y_0)^T ; U^\varepsilon(0)),
V^\varepsilon(s,\omega, (X_0, Y_0)^T ; U^\varepsilon(0)), \theta_s^\varepsilon\omega)ds\\
&=Y_0
+\frac{1}{\varepsilon}\int_{-\infty}^0e^{\frac{-Bs}{\varepsilon}}
\Delta G(U^\varepsilon(s, \omega, (X_0, Y_0)^T ; (\widetilde{X}_0-X_0)),\\
&~~~~~~~~~~~~~~~~~~~~~~~~~~~~~~~V^\varepsilon(s, \omega, (X_0,
Y_0)^T ; (\widetilde{X}_0-X_0)), \theta_s^\varepsilon\omega)ds,
\end{array}
$$
which just is $l^\varepsilon(\zeta,(X_0, Y_0)^T,\omega)$ if we take
$\widetilde{X}_0$ as $\zeta$ in $H_s$. Then from Lemma 4.1, Lemma
4.2, (\ref{Foliation}) and (\ref{Foliation-l}), we know that
$$
\mathcal{W}_\beta^\varepsilon((X_0, Y_0)^T, \omega)=\{(\zeta,
l^\varepsilon(\zeta,(X_0, Y_0)^T,\omega))^T| \; \zeta\in H_s\}.
$$
\par
Furthermore, for any $\zeta$ and $\widetilde{\zeta}$ in $H_s$, using
Lemma 4.3,
$$
\begin{array}{ll}
&\|l^\varepsilon(\zeta,(X_0,
Y_0)^T,\omega)-l^\varepsilon(\widetilde{\zeta},(X_0,
Y_0)^T,\omega)\|_{H_f}\\
=& \|V^\varepsilon(t,\omega, (X_0, Y_0)^T ;
\zeta-X_0)-V^\varepsilon(t,\omega, (X_0, Y_0)^T ;
\widetilde{\zeta}-X_0)\|_{H_f}|_{t=0}\\
\leq &\|V^\varepsilon(\cdot,\omega, (X_0, Y_0)^T ;
\zeta-X_0)-V^\varepsilon(\cdot,\omega, (X_0, Y_0)^T ; \widetilde{\zeta}-X_0)\|_{C_\beta^{f,-}}\\
\leq&
\frac{K}{\gamma_f+\varepsilon\beta}\|\Psi^\varepsilon(\cdot,\omega,
(X_0, Y_0)^T ; \zeta-X_0)-\Psi^\varepsilon(\cdot,\omega, (X_0,
Y_0)^T ;
\widetilde{\zeta}-X_0)\|_{C_\beta^{-}}\\
\leq&
\frac{K}{\gamma_f+\varepsilon\beta}\cdot\frac{1}{1-\rho(\gamma_s,\gamma_f,K,\varepsilon)}\|\zeta-\widetilde{\zeta}\|_{H_s}.
\end{array}
$$
\par
(ii) \emph{To prove the dynamical orbits are exponentially
approaching each other in backward time only if they start from the
same fiber.}
\par
From Lemma 4.1, using the same argument of the proof of Lemma 4.2,
we easily got
\begin{equation}\label{Theorem 3.1-1}
\begin{array}{lll}
\|\Psi^\varepsilon(\cdot)\|_{C_\beta^-}&=&\|U^\varepsilon(\cdot)\|_{C_\beta^{s,-}}
+\|V^\varepsilon(\cdot)\|_{C_\beta^{f,-}}\\
&=&\|e^{At}U^\varepsilon(0)\|_{C_\beta^{s,-}}+\|
\int_0^te^{A(t-s)}\Delta F(U^\varepsilon(s), V^\varepsilon(s),
\theta_s^\varepsilon\omega)ds\|_{C_\beta^{s,-}}\\
&&+\|
\frac{1}{\varepsilon}\int_{-\infty}^te^{\frac{B(t-s)}{\varepsilon}}
\Delta
G(U^\varepsilon(s), V^\varepsilon(s), \theta_s^\varepsilon\omega)ds\|_{C_\beta^{f,-}}\\
&\leq & \|U^\varepsilon(0)\|_{H_s}
+\frac{K}{-\beta-\gamma_s}\|\Psi^\varepsilon(\cdot)\|_{C_\beta^-}
+\frac{K}{\gamma_f+\varepsilon\beta}\|\Psi^\varepsilon(\cdot)\|_{C_\beta^-}\\
&\leq & \|U^\varepsilon(0)\|_{H_s}
+\rho(\gamma_s,\gamma_f,K,\varepsilon)\|\Psi^\varepsilon(\cdot)\|_{C_\beta^-},
\end{array}
\end{equation}
where $\Psi^\varepsilon$ is defined as (\ref{Equation-Psi}). Notice
that the Hypothesis H3 holds, $\beta=\frac{-\gamma_s}{2}$, and that
$\rho(\gamma_s,\gamma_f,K,\varepsilon)\longrightarrow
\frac{K}{-\beta-\gamma_s}+\frac{K}{\gamma_f}$ as
$\varepsilon\searrow 0$. Therefore, there exists a sufficiently
small positive constant $\varepsilon_0$ such that for
$\varepsilon\in (0, \varepsilon_0)$,
$\rho(\gamma_s,\gamma_f,K,\varepsilon)<1$. Then it follows from
(\ref{Theorem 3.1-1}) that
$$
\|\Psi^\varepsilon(\cdot)\|_{C_\beta^-}\leq
\frac{1}{1-\rho(\gamma_s,\gamma_f,K,\varepsilon)}
\|U^\varepsilon(0)\|_{H_s},
$$
which implies that
\begin{equation}\label{Theorem 3.2-2}
\|\Phi^\varepsilon(t,\omega,
(\widetilde{X}_0,\widetilde{Y}_0)^T)-\Phi^\varepsilon(t,\omega,
(X_0,Y_0)^T)\|_{H_s\times H_f}\leq \frac{e^{\beta
t}}{1-\rho(\gamma_s,\gamma_f,K,\varepsilon)}\cdot
\|U^\varepsilon(0)\|_{H_s}, \quad \forall \quad t\leq 0.
\end{equation}
\par
For any two points $(\widetilde{X}_0^1, \widetilde{Y}_0^1)^T$ and
$(\widetilde{X}_0^2, \widetilde{Y}_0^2)^T$ in the same fiber
$\mathcal{W}_\beta^\varepsilon((X_0, Y_0)^T, \omega)$, from
(\ref{Theorem 3.2-2}), we have
$$
\|\Phi^\varepsilon(t,\omega, (\widetilde{X}_0^1,
\widetilde{Y}_0^1)^T)-\Phi^\varepsilon(t,\omega,
(X_0,Y_0)^T)\|_{H_s\times H_f}\leq \frac{e^{\beta
t}}{1-\rho(\gamma_s,\gamma_f,K,\varepsilon)}\cdot
\|U^\varepsilon(0)\|_{H_s}, \quad \forall \quad t\leq 0,
$$
and
$$
\|\Phi^\varepsilon(t,\omega, (\widetilde{X}_0^2,
\widetilde{Y}_0^2)^T)-\Phi^\varepsilon(t,\omega,
(X_0,Y_0)^T)\|_{H_s\times H_f}\leq \frac{e^{\beta
t}}{1-\rho(\gamma_s,\gamma_f,K,\varepsilon)}\cdot
\|U^\varepsilon(0)\|_{H_s}, \quad \forall \quad t\leq 0,
$$
which immediately implies (\ref{Foliation-attracting}) holds.
\par
(iii) \emph{To prove the fiber is invariant.}
\par
To see this, taking a fiber
$\mathcal{W}_\beta^\varepsilon((X_0, Y_0)^T, \omega)$, we will show
that the time $\tau$-map $\Phi^\varepsilon(\tau,\omega,\cdot)$ maps
it into the fiber
$\mathcal{W}_\beta^\varepsilon(\Phi^\varepsilon(\tau, \omega, (X_0,
Y_0)^T), \theta_\tau\omega)$. Let $(\widetilde{X}_0,
\widetilde{Y}_0)^T\in \mathcal{W}_\beta^\varepsilon((X_0, Y_0)^T,
\omega)$. Then $\Phi^\varepsilon(\cdot, \omega, (\widetilde{X}_0,
\widetilde{Y}_0)^T)-\Phi^\varepsilon(\cdot, \omega, (X_0, Y_0)^T)\in
C_{\beta}^-$, which implies that
$$\Phi^\varepsilon(\cdot+\tau,
\omega, (\widetilde{X}_0,
\widetilde{Y}_0)^T)-\Phi^\varepsilon(\cdot+\tau, \omega, (X_0,
Y_0)^T)\in C_{\beta}^-.
$$
Thus by using the cocycle property
$$
\begin{array}{l}
\Phi^\varepsilon(\cdot +\tau, \omega, (\widetilde{X}_0,
\widetilde{Y}_0)^T)=\Phi^\varepsilon(\cdot, \theta_\tau\omega,
\Phi^\varepsilon(\tau, \omega, (\widetilde{X}_0,
\widetilde{Y}_0)^T)),\\
\Phi^\varepsilon(\cdot +\tau, \omega, (X_0,
Y_0)^T)=\Phi^\varepsilon(\cdot, \theta_\tau\omega,
\Phi^\varepsilon(\tau, \omega, (X_0, Y_0)^T)).
\end{array}
$$
Then we have $\Phi^\varepsilon(\tau, \omega, (X_0, Y_0)^T)\in
\mathcal{W}_\beta^\varepsilon(\Phi^\varepsilon(\tau, \omega, (X_0,
Y_0)^T), \theta_\tau\omega)$. The proof is completed.
\hfill{$\blacksquare$}
\par
{\bf Remark 4.1}\quad {\it From Theorem 4.1, \cite{WR-2013-JMAA} and
\cite{FLD}, we know that the invariance of the slow foliation means
the dynamical system maps a fiber to another fiber, while the
invariance of the slow manifold means the dynamical system preserve
the dynamical orbits starting from the slow manifold still in the
slow manifold.}
\par
\par
{\bf Remark 4.2}\quad {\it For the negative time, from Theorem 4.1,
we know that the slow foliation describes the dynamics of the system
(\ref{RE-s})-(\ref{RE-f}) in which the different dynamical orbits
are exponential closed only if they starting from a same fiber. For
the positive time, from \cite{WR-2013-JMAA} and \cite{FLD}, we know
that slow manifold describes the dynamics of the system
(\ref{RE-s})-(\ref{RE-f}), which could exponentially attract other
dynamical orbits. Therefore, the slow foliation and slow manifold
are from different view of points to describe the dynamics of the
slow-fast stochastic system. }
\par
{\bf Theorem 4.2 (Geometric properties of the slow foliation)}\quad
{\it Assume that the Hypothesis H1-H3 hold. Take $\beta$ as the
positive real number $\frac{-\gamma_s}{2}$. Let
$\mathcal{M}^\varepsilon(\omega)$
 and $\mathcal{W}_\beta^\varepsilon((X_0, Y_0)^T, \omega)$ be the slow manifold and a fiber
of the slow foliation for the slow-fast random system
(\ref{RE-s})-(\ref{RE-f}), respectively, which are well defined as
(\ref{manifold}) and (\ref{Foliation-mapping}).
Put
$$
\overset{\text{m}}{\mathcal{W}_\beta^\varepsilon}(\omega):=\{\mathcal{W}_\beta^\varepsilon((X_0,
Y_0)^T, \omega)\;|\; Y_0-h^\varepsilon(X_0, \omega):=m \in H_f,\;
(X_0, Y_0)^T\in H_s\times H_f\},
$$
where $h^\varepsilon(X_0, \omega)$ is defined as (\ref{manifold-h}).
%
%
Then the fiber
$\overset{\text{m}}{\mathcal{W}_\beta^\varepsilon}(\omega)$
parallels the fiber
$\overset{\text{n}}{\mathcal{W}^\varepsilon_\beta}(\omega)$ for any
$m, n\in H_f$ and $m\neq n$. Especially, if $m=0$,
$\overset{\text{m}}{\mathcal{W}^\varepsilon_\beta}(\omega)$ is just
the slow manifold. Thus, the geometry constructer of the slow
foliation is clear: every fiber of the slow foliation parallels each
other, and the slow manifold is one fiber of the slow foliation.
\par
Moreover, we have that
\par
(i) when the arbitrary given point $(X_0, Y_0)^T$ of
the slow foliation is in the slow manifold
$\mathcal{M}^\varepsilon(\omega)$, the fiber
$\mathcal{W}_\beta^\varepsilon((X_0, Y_0)^T, \omega)$ is just the
slow manifold $\mathcal{M}^\varepsilon(\omega)$;
\par
(ii) when the arbitrary given point $(X_0, Y_0)^T$ of the slow
foliation is not in the slow manifold
$\mathcal{M}^\varepsilon(\omega)$, the fiber
$\mathcal{W}_\beta^\varepsilon((X_0, Y_0)^T, \omega)$ parallels the
slow manifold $\mathcal{M}^\varepsilon(\omega)$. }
\par
{\bf Proof.}\quad From (\ref{Lemma 3.1-0}), for any
$(\widetilde{X}_0,\widetilde{Y}_0)^T\in H_s\times H_f$, we have
$$
\widetilde{Y}_0-Y_0 =
\frac{1}{\varepsilon}\int_{-\infty}^0e^{\frac{-Bs}{\varepsilon}}\Delta
G(U^\varepsilon(s), V^\varepsilon(s), \theta_s^\varepsilon\omega)ds,
$$
which implies that
\begin{equation}\label{Theorem 4.3-1-1}
\begin{array}{ll}
&\widetilde{Y}_0-
\frac{1}{\varepsilon}\int_{-\infty}^0e^{-\frac{Bs}{\varepsilon}}G(X^\varepsilon(s,
\omega; \widetilde{X}_0),  Y^\varepsilon(s, \omega;
\widetilde{X}_0), \theta_s^\varepsilon\omega)ds\\
=&Y_0 -
\frac{1}{\varepsilon}\int_{-\infty}^0e^{-\frac{Bs}{\varepsilon}}G(X^\varepsilon(s,
\omega; X_0),  Y^\varepsilon(s, \omega; X_0),
\theta_s^\varepsilon\omega)ds.
\end{array}
\end{equation}
In other words,
\begin{equation}\label{Theorem 4.3-1-2}
\widetilde{Y}_0-h^\varepsilon(\widetilde{X}_0,\omega)=Y_0-h^\varepsilon(X_0,\omega),
\end{equation}
where $h^\varepsilon(\cdot, \omega)$ is defined as
(\ref{manifold-h}).
\par
For arbitrary given point $(X_0, Y_0)^T$ of the slow foliation,
there exists $m\in H_f$ such that
$$
Y_0-h^\varepsilon(X_0,\omega)=m.
$$
\par
If $m=0$, then $(X_0, Y_0)^T$ is in the slow manifold
$\mathcal{M}^\varepsilon(\omega)$, which yields from (\ref{Theorem
4.3-1-2}) that
$$
\widetilde{Y}_0-h^\varepsilon(\widetilde{X}_0,\omega)=0, \quad
\hbox{for any}\quad  \widetilde{X}_0\in H_s.
$$
Thus,
$\overset{\text{0}}{\mathcal{W}^\varepsilon_\beta}(\omega)=\mathcal{M}^\varepsilon(\omega)$.
\par
If $m\neq 0$, then $(X_0, Y_0)^T$ is not in the slow manifold
$\mathcal{M}^\varepsilon(\omega)$. Then it immediately follows from
(\ref{Theorem 4.3-1-2}) that
$$
\widetilde{Y}_0-h^\varepsilon(\widetilde{X}_0,\omega)=m\neq0,\quad
\hbox{for any}\quad  \widetilde{X}_0\in H_s.
$$
Thus $(\widetilde{X}_0, \widetilde{Y}_0)^T $ is in the curve
$\overset{\text{m}}{\mathcal{W}_\beta^\varepsilon}(\omega)$ that
parallels the slow manifold
$\mathcal{M}^\varepsilon(\omega)=\overset{\text{0}}{\mathcal{W}^\varepsilon_\beta}(\omega)$.
Furthermore, for $m, n\in H_f$ and $m\neq n$, the
$\overset{\text{m}}{\mathcal{W}_\beta^\varepsilon}(\omega)$
parallels
$\overset{\text{n}}{\mathcal{W}^\varepsilon_\beta}(\omega)$. The
proof is completed.$\blacksquare$
\par

\renewcommand{\theequation}{\thesection.\arabic{equation}}
\setcounter{equation}{0}

\section{Critical foliation}

\quad\quad In this section, we will study the limiting case of the
slow foliation for the slow-fast random system
(\ref{RE-s})-(\ref{RE-f}) as the singular perturbation parameter
$\varepsilon$ tends to zero. Also, we delicately construct the
approximation of slow foliation for sufficiently small $\varepsilon$
in distribution.
\par
Taking the time scaling $t\to \varepsilon
t$ for the system (\ref{RE-s})-(\ref{RE-f}), we have
\begin{equation}\label{RE-s-scale}
\frac{dX^\varepsilon}{dt}=\varepsilon AX^\varepsilon+\varepsilon
F(X^\varepsilon, Y^\varepsilon,\theta_{\varepsilon
t}^\varepsilon\omega),
\end{equation}
\begin{equation}\label{RE-f-scale}
\frac{dY^\varepsilon}{dt}=BY^\varepsilon+G(X^\varepsilon,
Y^\varepsilon,\theta_{\varepsilon t}^\varepsilon\omega),~~
\end{equation}
where
$$
\begin{array}{l}
F(X^\varepsilon, Y^\varepsilon,\theta_{\varepsilon
t}^\varepsilon\omega):=f(X^\varepsilon+\delta(\theta_{\varepsilon
t}^1\omega_1),
Y^\varepsilon+\eta^{\frac{1}{\varepsilon}}(\theta_{\varepsilon
t}^2\omega_2)),\\
G(X^\varepsilon, Y^\varepsilon,\theta_{\varepsilon
t}^\varepsilon\omega):=g(X^\varepsilon+\delta(\theta_{\varepsilon
t}^1\omega_1),
Y^\varepsilon+\eta^{\frac{1}{\varepsilon}}(\theta_{\varepsilon
t}^2\omega_2)).
\end{array}
$$
\par
Noticing that Lemma 2.1, we replace
$\eta^{\frac{1}{\varepsilon}}(\theta_{\varepsilon t}^2\omega_2)$ by
$\xi(\theta_t^2\omega_2)$ in (\ref{RE-s-scale})-(\ref{RE-f-scale})
to get a new random evolutionary system
\begin{equation}\label{RE-s-distribution}
\frac{d\breve{X}^\varepsilon}{dt}=\varepsilon
A\breve{X}^\varepsilon+\varepsilon F(\breve{X}^\varepsilon,
\breve{Y}^\varepsilon,\overline{\theta_t^\varepsilon\omega}),
\end{equation}
\begin{equation}\label{RE-f-distribution}
\frac{d\breve{Y}^\varepsilon}{dt}=B\breve{Y}^\varepsilon+G(\breve{X}^\varepsilon,
\breve{Y}^\varepsilon,\overline{\theta_t^\varepsilon\omega}),~~
\end{equation}
where
$$
\begin{array}{l}
F(\breve{X}^\varepsilon,
\breve{Y}^\varepsilon,\overline{\theta_t^\varepsilon\omega}):=f(\breve{X}^\varepsilon+\delta(\theta_{\varepsilon
t}^1\omega_1),
\breve{Y}^\varepsilon+\xi(\theta_t^2\omega_2)),\\
G(\breve{X}^\varepsilon,
\breve{Y}^\varepsilon,\overline{\theta_t^\varepsilon\omega}):=
g(\breve{X}^\varepsilon+\delta(\theta_{\varepsilon t}^1\omega_1),
\breve{Y}^\varepsilon+\xi(\theta_t^2\omega_2)),
\end{array}
$$
with the initial condition $(\breve{X}^\varepsilon(0),
\breve{X}^\varepsilon(0))^T=(X_0, Y_0)^T$, whose solution is denoted
by
$$
\breve{\Phi}^\varepsilon(t,\omega, (X_0,
Y_0)^T)=(\breve{X}^\varepsilon(t, \omega, (X_0, Y_0)^T),
\breve{Y}^\varepsilon(t, \omega, (X_0, Y_0)^T))^T.
$$
Then the distribution of the solution
$\breve{\Phi}^\varepsilon(t,\omega, (X_0, Y_0)^T)$ of the system
(\ref{RE-s-distribution})-(\ref{RE-f-distribution}) coincides with
that of (\ref{RE-s})-(\ref{RE-f}) (also see \cite{SS}).
\par
Put
$$
\begin{array}{ll}
&\breve{\Phi}^\varepsilon(t,\omega,
(\widetilde{X}_0,\widetilde{Y}_0)^T)-\breve{\Phi}^\varepsilon(t,\omega,
(X_0,Y_0)^T)\\
= &(\breve{X}^\varepsilon(t,\omega,
(\widetilde{X}_0,\widetilde{Y}_0)^T)-\breve{X}^\varepsilon(t,\omega,
(X_0,Y_0)^T), \breve{Y}^\varepsilon(t,\omega,
(\widetilde{X}_0,\widetilde{Y}_0)^T)-\breve{Y}^\varepsilon(t,\omega,
(X_0,Y_0)^T))^T\\
=&:(\breve{U}^\varepsilon(t), \breve{V}^\varepsilon(t))^T.
\end{array}
$$
Then $\breve{U}^\varepsilon(t)=\breve{U}^\varepsilon(t,\omega, (X_0,
Y_0)^T ; (\widetilde{X}_0-X_0))$ and
$\breve{V}^\varepsilon(t)=\breve{V}^\varepsilon(t, \omega, (X_0,
Y_0)^T ; (\widetilde{X}_0-X_0))$ satisfies
\begin{equation}\label{Foliation-UV-distribution}
\left(
\begin{array}{l}
\breve{U}^\varepsilon(t)\\
\breve{V}^\varepsilon(t)
\end{array}
\right) = \left(
\begin{array}{c}
e^{\varepsilon A
t}\breve{U}^\varepsilon(0)+\varepsilon\int_0^te^{\varepsilon
A(t-s)}\Delta F(\breve{U}^\varepsilon(s),\breve{V}^\varepsilon(s),
\overline{\theta_s^\varepsilon\omega})ds\\
\int_{-\infty}^t e^{B(t-s)}\Delta
G(\breve{U}^\varepsilon(s),\breve{V}^\varepsilon(s),
\overline{\theta_s^\varepsilon\omega})ds
\end{array}
\right)
\end{equation}
with $\breve{U}^\varepsilon(0)=\widetilde{X}_0-X_0$. Here
$$
\begin{array}{ll}
\Delta F(\breve{U}^\varepsilon(s), \breve{V}^\varepsilon(s),
\overline{\theta_s^\varepsilon\omega}) =&
F(\breve{U}^\varepsilon(s)+\breve{X}^\varepsilon(s,\omega,
(X_0,Y_0)^T),
\breve{V}^\varepsilon(s)+\breve{Y}^\varepsilon(s,\omega,
(X_0,Y_0)^T),
\overline{\theta_s^\varepsilon\omega})\\
&-F(\breve{X}^\varepsilon(s,\omega, (X_0,Y_0)^T),
\breve{Y}^\varepsilon(s,\omega, (X_0,Y_0)^T),
\overline{\theta_s^\varepsilon\omega}),
\end{array}
$$
and
$$
\begin{array}{ll}
\Delta G(\breve{U}^\varepsilon(s), \breve{V}^\varepsilon(s),
\overline{\theta_s^\varepsilon\omega}) =&
G(\breve{U}^\varepsilon(s)+\breve{X}^\varepsilon(s,\omega,
(X_0,Y_0)^T),
\breve{V}^\varepsilon(s)+\breve{Y}^\varepsilon(s,\omega,
(X_0,Y_0)^T),
\overline{\theta_s^\varepsilon\omega})\\
&-G(\breve{X}^\varepsilon(s,\omega, (X_0,Y_0)^T),
\breve{Y}^\varepsilon(s,\omega, (X_0,Y_0)^T),
\overline{\theta_s^\varepsilon\omega}).
\end{array}
$$
\par
For every $\zeta\in H_s$, we define
\begin{equation}\label{Foliation-l-distribution}
\begin{array}{ll}
& \breve{l}^\varepsilon(\zeta,(X_0, Y_0)^T,\omega)\\
:=&Y_0 +\int_{-\infty}^0e^{-Bs} \Delta G(\breve{U}^\varepsilon(s,
\omega, (X_0, Y_0)^T ;
(\zeta-X_0)),\\
&~~~~~~~~~~~~~~~~~~~~~~~~~ \breve{V}^\varepsilon(s, \omega, (X_0,
Y_0)^T ; (\zeta-X_0)), \overline{\theta_s^\varepsilon\omega})ds.
\end{array}
\end{equation}
\par
Using the same arguments as in Section 4, we can obtain the slow
foliation of (\ref{RE-s-distribution})-(\ref{RE-f-distribution}) as
follows.
\par
{\bf Lemma 5.1}\quad {\it Assume that the Hypothesis H1-H3 hold.
Take $\beta$ as the positive real number $\frac{-\gamma_s}{2}$. Then
there exists a sufficiently small positive constant $\varepsilon_0$
such that for $\varepsilon\in (0, \varepsilon_0)$, the foliation of
the slow-fast random system
(\ref{RE-s-distribution})-(\ref{RE-f-distribution}) exists, whose
one fiber is given by
$$
\mathcal{\breve{W}}_\beta^\varepsilon((X_0, Y_0)^T,
\omega)=\{(\zeta, \breve{l}^\varepsilon(\zeta, (X_0, Y_0)^T,
\omega))^T| \; \zeta\in H_s\},
$$
where $(X_0, Y_0)^T\in H_s\times H_f$, the function
$\breve{l}^\varepsilon(\zeta,(X_0, Y_0)^T, \omega)$ is defined as
(\ref{Foliation-l-distribution}).

\par}

Furthermore, we obtain the relationship of foliation between of the
system (\ref{RE-s})-(\ref{RE-f}) and the system
(\ref{RE-s-distribution})-(\ref{RE-f-distribution}) as follows.
\par
{\bf Lemma 5.2}\quad {\it Assume that the Hypothesis H1-H3 hold.
Take $\beta$ as the positive real number $\frac{-\gamma_s}{2}$. Then
there exists a sufficiently small positive constant $\varepsilon_0$
such that for $\varepsilon\in (0, \varepsilon_0)$, the foliation of
the system (\ref{RE-s})-(\ref{RE-f}) is the same as that of the
system (\ref{RE-s-distribution})-(\ref{RE-f-distribution}) in
distribution, that is, for every $\zeta\in H_s$,
\begin{equation}\label{Lemma 4.2-0}
l^\varepsilon(\zeta,(X_0, Y_0)^T, \omega)\overset{\text{d}}{=}
\breve{l}^\varepsilon(\zeta,(X_0, Y_0)^T, \omega),
\end{equation}
where ``$\overset{\text{d}}{=}$" denotes equivalence in
distribution. }
\par
{\bf Proof.}\quad For $(\ref{Foliation-l})$, taking the time scaling
$s\to \varepsilon s$, and noticing that the solution of the system
(\ref{RE-s})-(\ref{RE-f}) has the same distribution as the solution
of the system (\ref{RE-s-distribution})-(\ref{RE-f-distribution}),
we know that for every $\zeta\in H_s$,
$$
\begin{array}{ll}
& l^\varepsilon(\zeta,(X_0, Y_0)^T,\omega)\\
=&Y_0
+\frac{1}{\varepsilon}\int_{-\infty}^0e^{\frac{-Bs}{\varepsilon}}
\Delta G(U^\varepsilon(s, \omega, (X_0, Y_0)^T ; (\zeta-X_0)),
V^\varepsilon(s, \omega, (X_0, Y_0)^T ; (\zeta-X_0)),
\theta_s^\varepsilon\omega)ds\\
=&Y_0 +\int_{-\infty}^0e^{-Bs} \Delta G(U^\varepsilon(\varepsilon s,
\omega, (X_0, Y_0)^T ; (\zeta-X_0)), V^\varepsilon(\varepsilon s,
\omega, (X_0, Y_0)^T ; (\zeta-X_0)),
\theta_{\varepsilon s}^\varepsilon\omega)ds\\
=&Y_0 +\int_{-\infty}^0e^{-Bs}[ g(U^\varepsilon(\varepsilon s,
\omega, (X_0, Y_0)^T ; (\zeta-X_0))+X^\varepsilon(\varepsilon
s,\omega,
(X_0,Y_0)^T)+\delta(\theta_{\varepsilon s}^1\omega_1), \\
&~~~~~~~~~~~~~~~~~~~~~~V^\varepsilon(\varepsilon s, \omega, (X_0,
Y_0)^T ; (\zeta-X_0))+Y^\varepsilon(\varepsilon s,\omega,
(X_0,Y_0)^T)+\eta^{\frac{1}{\varepsilon}}(
\theta_{\varepsilon s}^2\omega_2))\\
&~~~~~~~~~~~~~~~~~-g(X^\varepsilon(\varepsilon s,\omega,
(X_0,Y_0)^T)+\delta(\theta_{\varepsilon s}^1\omega_1),
Y^\varepsilon(\varepsilon s,\omega,
(X_0,Y_0)^T)+\eta^{\frac{1}{\varepsilon}} (\theta_{\varepsilon s}^2\omega_2))]ds\\
\overset{\text{d}}{=} &Y_0 +\int_{-\infty}^0e^{-Bs}[
g(\breve{U}^\varepsilon(s, \omega, (X_0, Y_0)^T ;
(\zeta-X_0))+\breve{X}^\varepsilon(s,\omega,
(X_0,Y_0)^T)+\delta(\theta_{\varepsilon s}^1\omega_1), \\
&~~~~~~~~~~~~~~~~~~~~~~\breve{V}^\varepsilon(s, \omega, (X_0, Y_0)^T
; (\zeta-X_0))+\breve{Y}^\varepsilon(s,\omega, (X_0,Y_0)^T)+\xi(
\theta_s^2\omega_2))\\
&~~~~~~~~~~~~~~~~~-g(\breve{X}^\varepsilon(s,\omega,
(X_0,Y_0)^T)+\delta(\theta_{\varepsilon s}^1\omega_1),
\breve{Y}^\varepsilon(s,\omega,
(X_0,Y_0)^T)+\xi(\theta_s^2\omega_2))]ds\\
=&Y_0 +\int_{-\infty}^0e^{-Bs} \Delta G(\breve{U}^\varepsilon(s,
\omega, (X_0, Y_0)^T ; (\zeta-X_0)), \breve{V}^\varepsilon(s,
\omega, (X_0, Y_0)^T ;
(\zeta-X_0)), \overline{\theta_s^\varepsilon\omega})ds\\
=&\breve{l}^\varepsilon(\zeta,(X_0, Y_0)^T,\omega).
\end{array}
$$
This completes the proof. \hfill{$\blacksquare$}
\par
Consider a new random evolutionary system
\begin{equation}\label{RE-s-critical}
\frac{dX^0}{dt}=0,~~~~~~~~~~~~~~~~~~~~~~~~~~~
\end{equation}
\begin{equation}\label{RE-f-critical}
\frac{d Y^0}{dt}=BY^0+G(X^0, Y^0,\overline{\theta_t^0\omega}),
\end{equation}
where
$$
G(X^0, Y^0, \overline{\theta_t^0\omega}):=g(X^0+\delta(\omega_1),
Y^0+\xi(\theta_t^2\omega_2)),
$$
with the initial condition $(X^0(0), Y^0(0))^T=(X_0, Y_0)^T$.
Essentially the system (\ref{RE-s-critical})-(\ref{RE-f-critical})
is the system (\ref{RE-s})-(\ref{RE-f}) scaled by $\varepsilon t$
with then zero singular perturbation parameter (i.e., the system
(\ref{RE-s-distribution})-(\ref{RE-f-distribution}) with
$\varepsilon=0$).
\par
We denote the solution of the system
(\ref{RE-s-critical})-(\ref{RE-f-critical}) as follows
$$
\Phi^0(t,\omega, (X_0, Y_0)^T)=(X_0, Y^0(t, \omega, (X_0,
Y_0)^T))^T.
$$
And put
$$
\begin{array}{ll}
&\Phi^0(t,\omega,
(\widetilde{X}_0,\widetilde{Y}_0)^T)-\Phi^0(t,\omega,
(X_0,Y_0)^T)\\
= &(\widetilde{X}_0-X_0, Y^0(t,\omega,
(\widetilde{X}_0,\widetilde{Y}_0)^T)-Y^0(t,\omega,
(X_0,Y_0)^T))^T\\
:=&(U^0(t), V^0(t))^T.
\end{array}
$$
Then $U^0(t)=U^0(0)=\widetilde{X}_0-X_0$ and $V^0(t)=V^0(t, \omega,
(X_0, Y_0)^T ; (\widetilde{X}_0-X_0))$ satisfies
\begin{equation}\label{Foliation-UV-critical}
\left(
\begin{array}{c}
U^0(t)\\
V^0(t)
\end{array}
\right) = \left(
\begin{array}{c}
U^0(0)\\
\int_{-\infty}^t e^{B(t-s)}\Delta G(U^0(0),V^0(s),
\overline{\theta_s^0\omega})ds
\end{array}
\right),
\end{equation}
where
$$
\begin{array}{lll}
\Delta G(U^0(0), V^0(s), \overline{\theta_s^0\omega}) &=&
G(U^0(0)+X_0, V^0(s)+Y^0(s,\omega, (X_0,Y_0)^T),
\overline{\theta_s^0\omega})\\
&&-G(X_0, Y^0(s,\omega, (X_0,Y_0)^T),\overline{\theta_s^0\omega}).
\end{array}
$$
For every $\zeta\in H_s$, we define
\begin{equation}\label{Foliation-l-critical}
\begin{array}{ll}
& l^0(\zeta,(X_0, Y_0)^T,\omega)\\
:=&Y_0 +\int_{-\infty}^0e^{-Bs} \Delta G((\zeta-X_0), V^0(s, \omega,
(X_0, Y_0)^T ; (\zeta-X_0)),\overline{\theta_s^0\omega})ds.
\end{array}
\end{equation}
\par
\par
Again using the same arguments as in Section 4, we can obtain the
invariant foliation of (\ref{RE-s-critical})-(\ref{RE-f-critical})
as follows.
\par
{\bf Theorem 5.1 (Critical foliation)}\quad {\it Assume that the
Hypothesis H1-H3 hold. Take $\beta$ as the positive real number
$\frac{-\gamma_s}{2}$. Then the invariant foliation of the random
system (\ref{RE-s-critical})-(\ref{RE-f-critical}) exists.
\par
(i)\quad Its one fiber is the graph of a Lipschitz function. That is
$$
\mathcal{W}_\beta^0((X_0, Y_0)^T, \omega)=\{(\zeta, l^0(\zeta, (X_0,
Y_0)^T, \omega))^T| \; \zeta\in H_s\},
$$
where $(X_0, Y_0)^T\in H_s\times H_f$, the function $l^0(\zeta,(X_0,
Y_0)^T, \omega)$ is defined as (\ref{Foliation-l-critical}). In
addition, $l^0(\zeta,(X_0, Y_0)^T, \omega)$ is Lipschitz continuous
with respect to $\zeta$, whose Lipschitz constant $Lip l^0$
satisfies
$$
Lip l^0 \leq \frac{K}{\gamma_f+\beta-K}.
$$
\par
(ii)\quad The dynamical orbits of
(\ref{RE-s-critical})-(\ref{RE-f-critical}) are exponentially
approaching each other in backward time only if they start from the
same fiber. That is, for any two points $(\widetilde{X}_0^1,
\widetilde{Y}_0^1)^T$ and $(\widetilde{X}_0^2, \widetilde{Y}_0^2)^T$
in a same fiber $\mathcal{W}_\beta^0((X_0, Y_0)^T, \omega)$,
$$
\begin{array}{ll}
\|\Phi^0(t,\omega, (\widetilde{X}_0^1,
\widetilde{Y}_0^1)^T)-\Phi^0(t,\omega, (\widetilde{X}_0^2,
\widetilde{Y}_0^2)^T)\|_{H_s\times H_f} =O(e^{\beta t}),\quad
\forall\quad t\to -\infty.
\end{array}
$$
\par
(iii)\quad Its fiber is invariant, i.e.,
$$
\Phi^0(t, \omega, \mathcal{W}_\beta^0((X_0, Y_0)^T, \omega))\subset
\mathcal{W}_\beta^0(\Phi^0(t, \omega, (X_0, Y_0)^T),
\overline{\theta_t^0\omega}).
$$
}
\par
{\bf Remark 5.1}\quad {\it From the Hypothesis H1-H3, and
$\beta=\frac{-\gamma_s}{2}$, we easily know that $\gamma_f+\beta-K$
is a positive constant.}
\par
\par
{\bf Remark 5.2}\quad {\it As we will show, the slow foliation of
the system (\ref{RE-s})-(\ref{RE-f}) converges to the foliation of
the system (\ref{RE-s-critical})-(\ref{RE-f-critical}) in
distribution, as $\varepsilon$ tends to zero. We call the limiting
status of the slow foliation as the\emph{ critical foliation} for
the system (\ref{RE-s})-(\ref{RE-f}). }
\par
{\bf Theorem 5.2 (Convergence in distribution to critical
foliation)}\quad{\it Assume that the Hypothesis H1-H3 hold. Take
$\beta$ as the positive real number $\frac{-\gamma_s}{2}$. And
assume that the nonlinear function $f(x, y)$ is bounded in $H_s$,
that is, there exists a positive constant such that  $\|f(x,
y)\|_{H_s}\leq C$. Then the slow foliation converges to the critical
foliation of the system of the system (\ref{RE-s})-(\ref{RE-f}) in
distribution (i.e., the distribution of the slow foliation converges
to the distribution of the critical foliation), as $\varepsilon$
tends to zero. In other words,
\begin{equation}\label{Theorem 5.2-0}
l^\varepsilon(\zeta, (X_0, Y_0)^T, \omega)\overset{\text{d}}{=}
l^0(\zeta, (X_0, Y_0)^T, \omega)+O(\varepsilon),\quad \hbox{in}\quad
H_f\quad \hbox{as}\quad \varepsilon\to 0.
\end{equation}
}
\par
{\bf Proof.}\quad Noticing that Lemma 5.2, we only need to prove
\begin{equation}\label{Theorem 5.2-1}
\breve{l}^\varepsilon(\zeta, (X_0, Y_0)^T,
\omega)\overset{\text{d}}{\longrightarrow} l^0(\zeta, (X_0, Y_0)^T,
\omega),\quad \hbox{in}\quad H_f\quad \hbox{as}\quad \varepsilon\to
0,
\end{equation}
for each given $\zeta\in H_s$.
\par
From (\ref{Foliation-l-distribution}) and
(\ref{Foliation-l-critical}), we know that
\begin{equation}\label{Theorem 5.2-2}
\begin{array}{ll}
&\|\breve{l}^\varepsilon(\zeta, (X_0, Y_0)^T,
\omega)-l^0(\zeta,(X_0, Y_0)^T,\omega)\|_{H_f}\\
=& \|\breve{V}^\varepsilon(t, \omega, (X_0, Y_0)^T ;
(\zeta-X_0))-V^0(t, \omega, (X_0, Y_0)^T ;
(\zeta-X_0))\|_{H_f}|_{t=0}.
\end{array}
\end{equation}
\par
For $\breve{V}^\varepsilon(t)$ and $V^0(t)$ with $t\leq 0$, we have
\begin{equation}\label{Theorem 5.2-3}
\begin{array}{ll}
&\|\breve{V}^\varepsilon(t, \omega, (X_0, Y_0)^T ; (\zeta-X_0))-V^0(t, \omega, (X_0, Y_0)^T ; (\zeta-X_0))\|_{H_f}\\
=& \|\int_{-\infty}^t e^{B(t-s)}[\Delta
G(\breve{U}^\varepsilon(s),\breve{V}^\varepsilon(s),
\overline{\theta_s^\varepsilon\omega})-\Delta G(U^0(0),V^0(s), \overline{\theta_s^0\omega})]ds\|_{H_f}\\
=& \|\int_{-\infty}^t e^{B(t-s)}[
G(\breve{U}^\varepsilon(s)+\breve{X}^\varepsilon(s,\omega,
(X_0,Y_0)^T),
\breve{V}^\varepsilon(s)+\breve{Y}^\varepsilon(s,\omega,
(X_0,Y_0)^T),
\overline{\theta_s^\varepsilon\omega})\\
&~~~~~~~~~~~~~~-G(U^0(0)+X_0, V^0(s)+Y^0(s,\omega, (X_0,Y_0)^T),
\overline{\theta_s^0\omega})\\
&~~~~~~~~~~~~~~-G(\breve{X}^\varepsilon(s,\omega, (X_0,Y_0)^T),
\breve{Y}^\varepsilon(s,\omega, (X_0,Y_0)^T), \overline{\theta_s^\varepsilon\omega})\\
&~~~~~~~~~~~~~~+G(X_0, Y^0(s,\omega, (X_0,Y_0)^T),
\overline{\theta_s^0\omega})
]ds\|_{H_f}\\
\leq & K  \int_{-\infty}^t e^{B(t-s)}[
\|\breve{U}^\varepsilon(s)-U^0(s)\|_{H_s}+\|\breve{V}^\varepsilon(s)-V^0(s)\|_{H_f}\\
&~~~~~~~~~~~~~~+2\|\breve{X}^\varepsilon(s)-X_0\|_{H_s}+2\|\breve{Y}^\varepsilon(s)-Y^0(s)\|_{H_f}
]ds
\end{array}
\end{equation}
for sufficiently small $\varepsilon$.
\par
To obtain the estimates of (\ref{Theorem 5.2-3}), we need to
establish the a priori estimates of
$\|\breve{X}^\varepsilon(t)-X_0\|_{H_s}$,
$\|\breve{Y}^\varepsilon(t)-Y^0(t)\|_{H_f}$, and
$\|\breve{U}^\varepsilon(t)-U^0(t)\|_{H_s}$, respectively.
\par
\emph{Step (i):} \emph{To estimate
$\|\breve{X}^\varepsilon(t)-X_0\|_{H_s}$.}
\par
For the system (\ref{RE-s-distribution})-(\ref{RE-f-distribution}),
using the same argument of Lemma 4.1, we can write it in a integral
form
\begin{equation}\label{RE-s-distribution-integral}
\left(\begin{array}{l} \breve{X}^\varepsilon(t)\\
\breve{Y}^\varepsilon(t)
\end{array}
\right) =\left(\begin{array}{c}
e^{\varepsilon
At}X_0+\varepsilon\int_0^t e^{\varepsilon
A(t-s)}F(\breve{X}^\varepsilon(s),
\breve{Y}^\varepsilon(s),\overline{\theta_s^\varepsilon\omega})ds\\
\int_{-\infty}^t e^{B(t-s)}G(\breve{X}^\varepsilon(s),
\breve{Y}^\varepsilon(s),\overline{\theta_s^\varepsilon\omega})ds
\end{array}
\right).
\end{equation}
\par
Then for any $t\leq 0$,
\begin{equation}\label{Theorem 5.2-4}
\begin{array}{ll}
&\|\breve{X}^\varepsilon(t)-X_0\|_{H_s}\\
\leq & \|e^{\varepsilon At}X_0-X_0\|_{H_s}+\varepsilon\int_t^0
e^{\varepsilon
\cdot(-\gamma_s)\cdot(t-s)}\|F(\breve{X}^\varepsilon(s),
\breve{Y}^\varepsilon(s),\overline{\theta_s^\varepsilon\omega})\|_{H_s}ds\\
\leq & \|e^{\varepsilon At}X_0-X_0\|_{H_s}+\varepsilon\cdot C
\cdot\int_t^0 e^{\varepsilon \cdot(-\gamma_s)\cdot(t-s)}ds\\
\leq & \|\int_{\varepsilon t}^0
AX_0e^{A\tau}d\tau\|_{H_s}+\varepsilon\cdot C
\cdot\int_t^0 e^{\varepsilon\cdot(-\gamma_s)\cdot(t-s)}ds\\
\leq & \|\int_{\varepsilon t}^0
AX_0e^{A\tau}d\tau\|_{H_s}+\varepsilon\cdot C
\cdot\int_t^0 e^{\varepsilon\cdot(-\gamma_s)\cdot(t-s)}ds\\
\leq &
\|AX_0\|_{H_s}\cdot\frac{1}{-\gamma_s}[1-e^{-\gamma_s\varepsilon
t}]+ C\cdot\frac{1}{-\gamma_s}[1-e^{-\gamma_s\varepsilon t}]\\
\leq & C[1-e^{-\gamma_s\varepsilon t}].
\end{array}
\end{equation}
Here and hereafter, we use $C$ to denote various positive constant
independent of $\varepsilon$ and $t$.
\par
\emph{Step (ii): }\emph{To estimate
$\|\breve{Y}^\varepsilon(t)-Y^0(t)\|_{H_f}$.}
\par
For the system (\ref{RE-s-critical})-(\ref{RE-f-critical}), using
the same argument of Lemma 4.1, we also can write it in a integral
form
\begin{equation}\label{RE-s-critical-integral}
\left(\begin{array}{l} X^0(t)\\
Y^0(t)
\end{array}
\right) =\left(\begin{array}{c} X_0\\
\int_{-\infty}^t e^{B(t-s)}G(X_0, Y^0(s),
\overline{\theta_s^0\omega})ds
\end{array}
\right).
\end{equation}
\par
Then for any $t\leq 0$, using (\ref{RE-s-distribution-integral}) and
(\ref{Theorem 5.2-4}), we deduce that
\begin{equation}\label{Theorem 5.2-5}
\begin{array}{ll}
&\|\breve{Y}^\varepsilon(t)-Y^0(t)\|_{H_f}\\
=&\| \int_{-\infty}^t e^{B(t-s)}[G(\breve{X}^\varepsilon(s),
\breve{Y}^\varepsilon(s),\overline{\theta_s^\varepsilon\omega})-G(X_0,
Y^0(s),\overline{\theta_s^0\omega})]ds \|_{H_f}\\
\leq & K  \int_{-\infty}^t
e^{-\gamma_f\cdot(t-s)}[\|\breve{X}^\varepsilon(s)-X_0\|_{H_s}+\|\breve{Y}^\varepsilon(s)-
Y^0(s)\|_{H_f}]ds  \\
\leq & K  \int_{-\infty}^t
e^{-\gamma_f\cdot(t-s)}[C(1-e^{-\gamma_s\cdot\varepsilon
s})+\|\breve{Y}^\varepsilon(s)- Y^0(s)\|_{H_f}]ds
\end{array}
\end{equation}
for sufficiently small $\varepsilon$.
\par
Take a real number $\alpha$ satisfying
\begin{equation}\label{Theorem 5.2-6}
\alpha\in [-\frac{2\gamma_f^2}{-\gamma_s+2\gamma_f}, 0).
\end{equation}
Then combining with the Hypothesis H1-H3, we know that
\begin{equation}\label{Theorem 5.2-7}
\alpha+\gamma_f>0,\quad -\alpha-\varepsilon\gamma_s>0, \quad
\hbox{and}\quad 0<\frac{K}{\gamma_f+\alpha}<1.
\end{equation}
\par
It follows from (\ref{Theorem 5.2-5}) that
\begin{equation}\label{Theorem 5.2-8}
\begin{array}{ll}
&\|\breve{Y}^\varepsilon-Y^0\|_{C_\alpha^{f,-}}\\
=&\sup\limits_{t\leq 0}e^{-\alpha t}\|\breve{Y}^\varepsilon(t)-Y^0(t)\|_{H_f}\\
\leq &  KC \sup\limits_{t\leq 0}e^{-\alpha t}\int_{-\infty}^t
e^{-\gamma_f\cdot(t-s)}[(1-e^{-\gamma_s\cdot\varepsilon s})ds+K
\sup\limits_{t\leq 0}e^{-\alpha t}\int_{-\infty}^t
e^{-\gamma_f\cdot(t-s)}\|\breve{Y}^\varepsilon(s)-
Y^0(s)\|_{H_f}ds\\
\leq &  KC \sup\limits_{t\leq 0}[\frac{1}{\gamma_f}e^{-\alpha
t}-\frac{1}{\gamma_f-\varepsilon\gamma_s}e^{-\alpha
t-\varepsilon\gamma_s t}]+
\frac{K}{\gamma_f+\alpha}\cdot\|\breve{Y}^\varepsilon-
Y^0\|_{C_\alpha^{f,-}}.
\end{array}
\end{equation}
\par
Define a function
\begin{equation}\label{Theorem 5.2-9}
p(t):=\frac{1}{\gamma_f}e^{-\alpha
t}-\frac{1}{\gamma_f-\varepsilon\gamma_s}e^{-\alpha
t-\varepsilon\gamma_s t}, \quad\hbox{for}\quad t\leq 0.
\end{equation}
Then
$$
\begin{array}{ll}
p^\prime(t)&=e^{-\alpha
t}[\frac{-\alpha}{\gamma_f}-\frac{-\alpha-\varepsilon
\gamma_s}{\gamma_f-\varepsilon \gamma_2}e^{-\varepsilon \gamma_s
t}]\\
&\geq e^{-\alpha
t}[\frac{-\alpha}{\gamma_f}-\frac{-\alpha-\varepsilon
\gamma_s}{\gamma_f-\varepsilon \gamma_2}e^{0}]\\
&\longrightarrow 0, \quad \hbox{as}\quad \varepsilon \to 0,
\end{array}
$$
which implies there exists a sufficient small $\varepsilon_0$ such
that for $\varepsilon\in (0, \varepsilon_0)$, $p(t)$ is increasing
with respect to the variable $t$. Then we immediately get
\begin{equation}\label{Theorem 5.2-10}
p(t)\leq p(0)=\frac{1}{\gamma_f}-\frac{1}{\gamma_f-\varepsilon
\gamma_s},\quad \hbox{for} \quad t\leq 0.
\end{equation}
Thus, it follows from (\ref{Theorem 5.2-8})-(\ref{Theorem 5.2-10})
that
\begin{equation}\label{Theorem 5.2-11}
\|\breve{Y}^\varepsilon-Y^0\|_{C_\alpha^{f,-}}\leq
\frac{KC}{1-\frac{K}{\gamma_f+\alpha}}(\frac{1}{\gamma_f}-\frac{1}{\gamma_f-\varepsilon
\gamma_s}).
\end{equation}
Then
\begin{equation}\label{Theorem 5.2-12}
\|\breve{Y}^\varepsilon(t)-Y^0(t)\|_{H_f}\leq
\frac{KC}{1-\frac{K}{\gamma_f+\alpha}}(\frac{1}{\gamma_f}-\frac{1}{\gamma_f-\varepsilon
\gamma_s})\cdot e^{\alpha t},\quad \hbox{for} \quad t\leq 0,
\end{equation}
which is significative from (\ref{Theorem 5.2-6}) and (\ref{Theorem
5.2-7}).
\par
\emph{Step (iii):} \emph{To estimate
$\|\breve{U}^\varepsilon(t)-U^0(t)\|_{H_s}$.}
\par
It follows from (\ref{Foliation-UV-distribution}) and
(\ref{Foliation-UV-critical}) that
$$
\begin{array}{ll}
&\|\breve{U}^\varepsilon(t)-U^0(t)\|_{H_s}\\
\leq &\| e^{\varepsilon A
t}\breve{U}^\varepsilon(0)-U^0(0)\|_{H_s}+\|\int_0^te^{\varepsilon
A(t-s)}\Delta F(\breve{U}^\varepsilon(s),\breve{V}^\varepsilon(s),
\overline{\theta_s^\varepsilon\omega})ds\|_{H_s}.
\end{array}
$$
Using the same argument as (i), we can get
\begin{equation}\label{Theorem 5.2-13}
\|\breve{U}^\varepsilon(t)-U^0(t)\|_{H_s}\leq
C[1-e^{-\gamma_s\varepsilon t}].
\end{equation}
\par
\vspace{0.5cm}
\par
Now we go back to (\ref{Theorem 5.2-3}) to estimate
$\|\breve{V}^\varepsilon(t)-V^0(t)\|_{H_f}$.
\par
It follows from (\ref{Theorem 5.2-3}), (\ref{Theorem 5.2-4}),
(\ref{Theorem 5.2-12}) and (\ref{Theorem 5.2-13}) that for any
$t\leq 0$,
$$
\begin{array}{ll}
&\|\breve{V}^\varepsilon(t, \omega, (X_0, Y_0)^T ; (\zeta-X_0))-V^0(t, \omega, (X_0, Y_0)^T ; (\zeta-X_0))\|_{H_f}\\
\leq & K \int_{-\infty}^t e^{-\gamma_f\cdot(t-s)}[
3C[1-e^{-\gamma_s\varepsilon
s}]+2\frac{KC}{1-\frac{K}{\gamma_f+\alpha}}(\frac{1}{\gamma_f}-\frac{1}{\gamma_f-\varepsilon
\gamma_s})\cdot e^{\alpha
s}+\|\breve{V}^\varepsilon(s)-V^0(s)\|_{H_f} ]ds,
\end{array}
$$
which implies that
\begin{equation}\label{Theorem 5.2-14}
\begin{array}{ll}
&\|\breve{V}^\varepsilon-V^0\|_{C_\alpha^{f,-}}\\
\leq
&K\|\breve{V}^\varepsilon-V^0\|_{C_\alpha^{f,-}}\cdot\sup\limits_{t\leq
0}e^{-\alpha t} \int_{-\infty}^t e^{-\gamma_f\cdot(t-s)}e^{\alpha s}
ds\\
&+KC\cdot\sup\limits_{t\leq 0}e^{-\alpha t} \int_{-\infty}^t
e^{-\gamma_f\cdot(t-s)}[(1-e^{-\gamma_s\varepsilon
s})+(\frac{1}{\gamma_f}-\frac{1}{\gamma_f-\varepsilon
\gamma_s})\cdot e^{\alpha s}]ds\\
\leq
&\frac{K}{\gamma_f+\alpha}\cdot\|\breve{V}^\varepsilon-V^0\|_{C_\alpha^{f,-}}\\
&+KC\cdot \sup\limits_{t\leq
0}[\frac{1}{\gamma_f+\alpha}(\frac{1}{\gamma_f}-\frac{1}{\gamma_f-\varepsilon
\gamma_s})+\frac{1}{\gamma_f}e^{-\alpha
t}-\frac{1}{\gamma_f-\varepsilon \gamma_s} e^{-\alpha t-\varepsilon
\gamma_s t}].
\end{array}
\end{equation}
\par
Again using (\ref{Theorem 5.2-9}) and (\ref{Theorem 5.2-10}), then
it follows from (\ref{Theorem 5.2-14}) that
\begin{equation}\label{Theorem 5.2-15}
\|\breve{V}^\varepsilon-V^0\|_{C_\alpha^{f,-}} \leq
\frac{KC(1+\frac{1}{\gamma_f+\alpha})}{1-\frac{K}{\gamma_f+\alpha}}\cdot
(\frac{1}{\gamma_f}-\frac{1}{\gamma_f-\varepsilon \gamma_s}),
\end{equation}
which is also significative from (\ref{Theorem 5.2-6}) and
(\ref{Theorem 5.2-7}). Then we immediately have that
\begin{equation}\label{Theorem 5.2-16}
\|\breve{V}^\varepsilon(t)-V^0(t)\|_{H_f} \leq
\frac{KC(1+\frac{1}{\gamma_f+\alpha})}{1-\frac{K}{\gamma_f+\alpha}}\cdot
(\frac{1}{\gamma_f}-\frac{1}{\gamma_f-\varepsilon \gamma_s})\cdot
e^{\alpha t}, \quad\hbox{for}\quad t\leq 0.
\end{equation}
\par
Hence, it finally follows from (\ref{Theorem 5.2-2}) and
(\ref{Theorem 5.2-16}) that
\begin{equation}\label{Theorem 5.2-17}
\begin{array}{ll}
&\|\breve{l}^\varepsilon(\zeta, (X_0, Y_0)^T,
\omega)-l^0(\zeta,(X_0, Y_0)^T,\omega)\|_{H_f}\\
=& \|\breve{V}^\varepsilon(t, \omega, (X_0, Y_0)^T ; (\zeta-X_0))-V^0(t, \omega, (X_0, Y_0)^T ; (\zeta-X_0))\|_{H_f}|_{t=0}\\
\leq &
\frac{KC(1+\frac{1}{\gamma_f+\alpha})}{1-\frac{K}{\gamma_f+\alpha}}\cdot
(\frac{1}{\gamma_f}-\frac{1}{\gamma_f-\varepsilon \gamma_s})\\
\longrightarrow & 0,\quad \hbox{as}\quad \varepsilon \to 0,
\end{array}
\end{equation}
which implies (\ref{Theorem 5.2-1}) holds. This completes the proof.
\hfill{$\blacksquare$}
\par
{\bf Theorem 5.3 (Approximation of slow foliation)}\quad{\it Assume
that the Hypothesis H1-H3 hold. Take $\beta$ as the positive real
number $\frac{-\gamma_s}{2}$. And assume that the nonlinear function
$f(x, y)$ is bounded in $H_s$. Then for sufficiently small
$\varepsilon$, the slow foliation of the system
(\ref{RE-s})-(\ref{RE-f}) can be approximated in distribution as
\begin{equation}\label{Theorem 5.3-0}
\begin{array}{ll}
\mathcal{W}_\beta^\varepsilon((X_0, Y_0)^T, \omega)&=\{(\zeta,
l^\varepsilon(\zeta, (X_0, Y_0)^T, \omega))|\; \zeta\in
H_s\}\\
&\overset{\text{d}}{=}\{(\zeta, l^0(\zeta, (X_0, Y_0)^T,
\omega)+\varepsilon l^1(\zeta, (X_0, Y_0)^T,
\omega)+O(\varepsilon^2))|\; \zeta\in H_s\},
\end{array}
\end{equation}
where $l^0(\zeta, (X_0, Y_0)^T, \omega)$ is the critical foliation
as (\ref{Foliation-l-critical}), and $l^1(\zeta, (X_0, Y_0)^T,
\omega)$ is well defined as (\ref{foliation-l-approximation}). }
\par
{\bf Proof.}\quad From Lemma 5.2, it is only need to prove
\begin{equation}\label{Theorem 5.3-1}
\breve{l}^\varepsilon(\zeta, (X_0, Y_0)^T, \omega)=l^0(\zeta, (X_0,
Y_0)^T, \omega)+\varepsilon l^1(\zeta, (X_0, Y_0)^T,
\omega)+O(\varepsilon^2),\quad \hbox{in}\quad H_f.
\end{equation}
\par
For the system (\ref{RE-s-distribution})-(\ref{RE-f-distribution}),
we write
\begin{equation}\label{Theorem 5.3-2}
\begin{array}{l}
\breve{X}^\varepsilon(t)=\breve{X}^0(t)+\varepsilon
X^1(t)+O(\varepsilon^2),\\
\breve{X}^\varepsilon(0)=X_0,
\end{array}
\end{equation}
and
\begin{equation}\label{Theorem 5.3-3}
\begin{array}{l}
\breve{Y}^\varepsilon(t)=\breve{Y}^0(t)+\varepsilon
Y^1(t)+O(\varepsilon^2),\\
\breve{Y}^\varepsilon(0)=Y_0,
\end{array}
\end{equation}
where $\breve{X}^0(t)$, $\breve{Y}^0(t)$, $X^1(t)$ and $Y^1(t)$ will
be determined in the below. Also, notice that the relationship of
$(\breve{X}^\varepsilon(t), \breve{Y}^\varepsilon(t))^T$ and
$(\breve{U}^\varepsilon(t), \breve{V}^\varepsilon(t))^T$. We can
write
\begin{equation}\label{Theorem 5.3-4}
\begin{array}{l}
\breve{U}^\varepsilon(t)=\breve{U}^0(t)+\varepsilon
U^1(t)+O(\varepsilon^2),\\
\breve{U}^\varepsilon(0)=\zeta-X_0,
\end{array}
\end{equation}
and
\begin{equation}\label{Theorem 5.3-5}
\begin{array}{l}
\breve{V}^\varepsilon(t)=\breve{V}^0(t)+\varepsilon
V^1(t)+O(\varepsilon^2),\\
\breve{V}^\varepsilon(0)=\breve{l}(\zeta, (X_0, Y_0)^T, \omega)-Y_0.
\end{array}
\end{equation}

Expanding $F(\breve{X}^\varepsilon(t), \breve{Y}^\varepsilon(t),
\overline{\theta_t^\varepsilon\omega})$ at $\varepsilon=0$ by Taylor
formula, we infer that
\begin{equation}\label{Theorem 5.3-6}
\begin{array}{ll}
&F(\breve{X}^\varepsilon(t), \breve{Y}^\varepsilon(t),
\overline{\theta_t^\varepsilon\omega})\\
=& f(\breve{X}^\varepsilon(t)+\delta(\theta_{\varepsilon
t}^1\omega_1),
\breve{Y}^\varepsilon(t)+\xi(\theta_{t}^2\omega_2))\\
=& f(\breve{X}^0(t)+\delta(\omega_1),
\breve{Y}^0(t)+\xi(\theta_{t}^2\omega_2))\\
&+(\breve{X}^\varepsilon(t)-\breve{X}^0(t))f_x(\breve{X}^0(t)+\delta(\omega_1),
\breve{Y}^0(t)+\xi(\theta_{t}^2\omega_2))\\
&+(\breve{Y}^\varepsilon(t)-\breve{Y}^0(t))f_y(\breve{X}^0(t)+\delta(\omega_1),
\breve{Y}^0(t)+\xi(\theta_{t}^2\omega_2))+O(\varepsilon^2)\\
=& F(\breve{X}^0(t),
\breve{Y}^0(t), \overline{\theta_t^0\omega})\\
&+\varepsilon\cdot X^1(t)\cdot f_x(\breve{X}^0(t)+\delta(\omega_1),
\breve{Y}^0(t)+\xi(\theta_{t}^2\omega_2))\\
&+\varepsilon\cdot Y^1(t)\cdot f_y(\breve{X}^0(t)+\delta(\omega_1),
\breve{Y}^0(t)+\xi(\theta_{t}^2\omega_2))+O(\varepsilon^2),
\end{array}
\end{equation}
where $f_x(\cdot, \cdot)$ and $f_y(\cdot, \cdot)$ denote the partial
derivative of $f(x,y)$ with respect to the first variable $x$, and
the second variable $y$, respectively.
\par
Similarly, we get
\begin{equation}\label{Theorem 5.3-7}
\begin{array}{ll}
&G(\breve{X}^\varepsilon(t), \breve{Y}^\varepsilon(t),
\overline{\theta_t^\varepsilon\omega})\\
=& G(\breve{X}^0(t),
\breve{Y}^0(t), \overline{\theta_t^0\omega})\\
&+\varepsilon\cdot X^1(t)\cdot g_x(\breve{X}^0(t)+\delta(\omega_1),
\breve{Y}^0(t)+\xi(\theta_{t}^2\omega_2))\\
&+\varepsilon\cdot Y^1(t)\cdot g_y(\breve{X}^0(t)+\delta(\omega_1),
\breve{Y}^0(t)+\xi(\theta_{t}^2\omega_2))+O(\varepsilon^2),
\end{array}
\end{equation}
where $g_x(\cdot, \cdot)$ and $g_y(\cdot, \cdot)$ denote the partial
derivative of $g(x,y)$ with respect to the first variable $x$, and
the second variable $y$, respectively. We also have
\begin{equation}\label{Theorem 5.3-8}
\begin{array}{ll}
&F(\breve{U}^\varepsilon(t)+\breve{X}^\varepsilon(t),
\breve{V}^\varepsilon(t)+\breve{Y}^\varepsilon(t),
\overline{\theta_t^\varepsilon\omega})\\
=& F(\breve{U}^0(t)+\breve{X}^0(t),
\breve{V}^0(t)+\breve{Y}^0(t), \overline{\theta_t^0\omega})\\
&+\varepsilon\cdot (U^1(t)+X^1(t))\cdot
f_x(\breve{U}^0(t)+\breve{X}^0(t)+\delta(\omega_1),
\breve{V}^0(t)+\breve{Y}^0(t)+\xi(\theta_{t}^2\omega_2))\\
&+\varepsilon\cdot (V^1(t)+Y^1(t))\cdot
f_y(\breve{U}^0(t)+\breve{X}^0(t)+\delta(\omega_1),\breve{V}^0(t)+
\breve{Y}^0(t)+\xi(\theta_{t}^2\omega_2))\\
&+O(\varepsilon^2),
\end{array}
\end{equation}
and
\begin{equation}\label{Theorem 5.3-9}
\begin{array}{ll}
&G(\breve{U}^\varepsilon(t)+\breve{X}^\varepsilon(t),
\breve{V}^\varepsilon(t)+\breve{Y}^\varepsilon(t),
\overline{\theta_t^\varepsilon\omega})\\
=& G(\breve{U}^0(t)+\breve{X}^0(t),
\breve{V}^0(t)+\breve{Y}^0(t), \overline{\theta_t^0\omega})\\
&+\varepsilon\cdot (U^1(t)+X^1(t))\cdot
g_x(\breve{U}^0(t)+\breve{X}^0(t)+\delta(\omega_1),
\breve{V}^0(t)+\breve{Y}^0(t)+\xi(\theta_{t}^2\omega_2))\\
&+\varepsilon\cdot (V^1(t)+Y^1(t))\cdot
g_y(\breve{U}^0(t)+\breve{X}^0(t)+\delta(\omega_1),\breve{V}^0(t)+
\breve{Y}^0(t)+\xi(\theta_{t}^2\omega_2))\\
&+O(\varepsilon^2).
\end{array}
\end{equation}

Substituting (\ref{Theorem 5.3-2}) into (\ref{RE-s-distribution}),
and substituting (\ref{Theorem 5.3-3}) into
(\ref{RE-f-distribution}), then equating the terms with the same
power of $\varepsilon$, we deduce that
\begin{equation}\label{Theorem 5.3-10}
\begin{array}{l}
\breve{X}_t^0(t)=0,\\
\breve{Y}_t^0(t)=B\breve{Y}^0(t)+G(\breve{X}^0(t), \breve{Y}^0(t),
\overline{\theta_t^0\omega})
\end{array}
\end{equation}
and
\begin{equation}\label{Theorem 5.3-11}
\begin{array}{ll}
X_t^1(t))=&A\breve{X}^0(t)+F(\breve{X}^0(t), \breve{Y}^0(t),
\overline{\theta_t^0\omega}).\\
Y_t^1(t))=&BY^1(t)+X^1(t)\cdot g_x(\breve{X}^0(t)+\delta(\omega_1),
\breve{Y}^0(t)+\xi(\theta_t^2\omega_2))\\
&~~~~~~~~~~+Y^1(t)\cdot g_y(\breve{X}^0(t)+\delta(\omega_1),
\breve{Y}^0(t)+\xi(\theta_t^2\omega_2)).
\end{array}
\end{equation}
Similarly, we have
\begin{equation}\label{Theorem 5.3-12}
\begin{array}{ll}
\breve{U}_t^0(t)=&0,\\
\breve{V}_t^0(t)=&B\breve{V}^0(t)+[G(\breve{U}^0(t)+\breve{X}^0(t),
\breve{V}^0(t)+\breve{Y}^0(t),
\overline{\theta_t^0\omega})\\
&~~~~~~~~~~~~~-G(\breve{X}^0(t), \breve{Y}^0(t),
\overline{\theta_t^0\omega})]
\end{array}
\end{equation}
and
\begin{equation}\label{Theorem 5.3-13}
\begin{array}{ll}
U_t^1(t))=&A\breve{U}^0(t)+[F(\breve{U}^0(t)+\breve{X}^0(t),
\breve{V}^0(t)+\breve{Y}^0(t), \overline{\theta_t^0\omega})
-F(\breve{X}^0(t), \breve{Y}^0(t), \overline{\theta_t^0\omega})].\\
V_t^1(t))=&B V^1(t)+[(U^1(t)+X^1(t))\cdot
g_x(\breve{U}^0(t)+\breve{X}^0(t)+\delta(\omega_1),
\breve{V}^0(t)+\breve{Y}^0(t)+\xi(\theta_t^2\omega_2))\\
&~~~~~~~~~~+(V^1(t)+Y^1(t))\cdot
g_y(\breve{U}^0(t)+\breve{X}^0(t)+\delta(\omega_1),
\breve{V}^0(t)+\breve{Y}^0(t)+\xi(\theta_t^2\omega_2))\\
&~~~~~~~~~~-X^1(t)\cdot g_x(\breve{X}^0(t)+\delta(\omega_1),
\breve{Y}^0(t)+\xi(\theta_t^2\omega_2))\\
&~~~~~~~~~~-Y^1(t)\cdot g_y(\breve{X}^0(t)+\delta(\omega_1),
\breve{Y}^0(t)+\xi(\theta_t^2\omega_2))].
\end{array}
\end{equation}
\par
It immediately follows from (\ref{RE-s-critical}),
(\ref{Foliation-UV-critical}), (\ref{Theorem 5.3-10}) and
(\ref{Theorem 5.3-12}) that
\begin{equation}\label{Theorem 5.3-14}
\breve{X}^0(t)=X^0(t),\quad \breve{Y}^0(t)=Y^0(t),\quad
\breve{U}^0(t)=U^0(t), \quad \breve{V}^0(t)=V^0(t).
\end{equation}
In addition, using the contraction mapping principle as in Lemma
4.2, we can easily obtain the existence of $(X^1(t), Y^1(t))^T$ and
$(U^1(t), V^1(t))^T$. Here, for simplicity, we omit it. Then we can
define
\begin{equation}\label{foliation-l-approximation}
\begin{array}{ll}
l^1(\zeta, (X_0, Y_0)^T, \omega)=&\int_{-\infty}^0e^{-Bs}
[(U^1(t)+X^1(t))\cdot g_x(\zeta+\delta(\omega_1),
V^0(t)+Y^0(t)+\xi(\theta_t^2\omega_2))\\
&~~~~~~~~~~+(V^1(t)+Y^1(t))\cdot g_y(\zeta+\delta(\omega_1),
V^0(t)+Y^0(t)+\xi(\theta_t^2\omega_2))\\
&~~~~~~~~~~-X^1(t)\cdot g_x(X_0+\delta(\omega_1),
Y^0(t)+\xi(\theta_t^2\omega_2))\\
&~~~~~~~~~~-Y^1(t)\cdot g_y(X_0+\delta(\omega_1),
Y^0(t)+\xi(\theta_t^2\omega_2))]ds.
\end{array}
\end{equation}
\par
For every $\zeta\in H_s$, noticing that the initial condition
$X^0(0)=X_0$ and $U^0(0)=\zeta-X_0$, it follows from
(\ref{Foliation-l-distribution}), (\ref{Theorem 5.3-7}),
(\ref{Theorem 5.3-9}) and (\ref{Theorem 5.3-14}) that
\begin{equation}\label{Theorem 5.3-15}
\begin{array}{ll}
& \breve{l}^\varepsilon(\zeta,(X_0, Y_0)^T,\omega)\\
=&Y_0 +\int_{-\infty}^0e^{-Bs} [
G(\breve{U}^\varepsilon(s)+\breve{X}^\varepsilon(s),
\breve{V}^\varepsilon(s)+\breve{Y}^\varepsilon(s),\overline{\theta_s^\varepsilon\omega})\\
&~~~~~~~~~~~~~~~~~~~~~~~~~ -G(\breve{X}^\varepsilon(s),
\breve{Y}^\varepsilon(s),\overline{\theta_s^\varepsilon\omega})]ds\\
=&Y_0 +\int_{-\infty}^0e^{-Bs} \bigtriangleup G(U^0(s),
V^0(s),\overline{\theta_s^0\omega})ds\\
&+\varepsilon\int_{-\infty}^0e^{-Bs} [(U^1(t)+X^1(t))\cdot
g_x(U^0(t)+X^0(t)+\delta(\omega_1),
V^0(t)+Y^0(t)+\xi(\theta_t^2\omega_2))\\
&~~~~~~~~~~+(V^1(t)+Y^1(t))\cdot g_y(U^0(t)+X^0(t)+\delta(\omega_1),
V^0(t)+Y^0(t)+\xi(\theta_t^2\omega_2))\\
&~~~~~~~~~~-X^1(t)\cdot g_x(X^0(t)+\delta(\omega_1),
Y^0(t)+\xi(\theta_t^2\omega_2))\\
&~~~~~~~~~~-Y^1(t)\cdot g_y(X^0(t)+\delta(\omega_1),
Y^0(t)+\xi(\theta_t^2\omega_2))]ds\\
&+O(\varepsilon^2)\\
=&l^0(\zeta, (X_0, Y_0)^T, \omega)\\
&+\varepsilon\int_{-\infty}^0e^{-Bs} [(U^1(t)+X^1(t))\cdot
g_x(\zeta+\delta(\omega_1),
V^0(t)+Y^0(t)+\xi(\theta_t^2\omega_2))\\
&~~~~~~~~~~+(V^1(t)+Y^1(t))\cdot g_y(\zeta+\delta(\omega_1),
V^0(t)+Y^0(t)+\xi(\theta_t^2\omega_2))\\
&~~~~~~~~~~-X^1(t)\cdot g_x(X_0+\delta(\omega_1),
Y^0(t)+\xi(\theta_t^2\omega_2))\\
&~~~~~~~~~~-Y^1(t)\cdot g_y(X_0+\delta(\omega_1),
Y^0(t)+\xi(\theta_t^2\omega_2))]ds\\
&+O(\varepsilon^2),
\end{array}
\end{equation}
which immediately follows from (\ref{foliation-l-approximation})
that Theorem 4.3 holds. \hfill{$\blacksquare$}
\par
{\bf Remark 5.3}\quad  {\it In Section 4 and Section 5, the
conditions for the general theory of the slow foliation are only
sufficient condition, not the necessary condition. }
\par
{\bf Example 5.1}\quad  {\it Consider the following slow-fast
stochastic evolutionary system
\begin{equation}\label{Equation-s-example 5.1}
\frac{dx^\varepsilon}{dt}=x^\varepsilon+f(x^\varepsilon,
y^\varepsilon)+\sigma_1\dot{W_1},~~~~~~~~~~~\quad \hbox{in}\quad
H_s,
\end{equation}
\begin{equation}\label{Equation-f-example 5.1}
\frac{dy^\varepsilon}{dt}=\frac{1}{\varepsilon} \bigtriangleup
y^\varepsilon+\frac{1}{\varepsilon}g(x^\varepsilon,
y^\varepsilon)+\frac{\sigma_2}{\sqrt{\varepsilon}}\dot{W_2},\quad
\hbox{in}\quad H_f.
\end{equation}
The system may model certain biological processes, for instance, the
famous FitzHugh-Nagumo system, as a simplified version of the
Hodgkin-Huxley model, which describes mechanisms of a neural
excitability and excitation for macro-receptors.
\par
Let $A$ be $\hbox{Id}$ (the identity operator) in
$H_s=L^2([0,\pi])$. Then it is clear that $\|e^{At}x\|_{H_s}\leq
e^{t}\|x\|_{H_s}= e^{-\gamma_s\cdot t}\|x\|_{H_s}$ with
$\gamma_s=-1$. Let $B$ be $\bigtriangleup$ with domain
$\mathcal{D}=H^2([0,\pi])\bigcap H^1_0([0,\pi])$, whose eigenvalue
are $\lambda_k=-k^2$ with the corresponding eigenfunction $e_k=\sin
kx$ $(k=1,2,\cdots)$, generating a $C_0$-semigroup $\{e^{Bt}: t\geq
0\}$ on $H_f=L^2([0,\pi])$ satisfying $\|e^{Bt}y\|_{H_f}\leq
e^{-\gamma_f\cdot t}\|y\|_{H_f}$ with $\gamma_f=1$.
\par
Assume that nonlinear functions $f: H_s\times H_f\longrightarrow H_s
$ and $g: H_s\times H_f\longrightarrow H_f $, which are $C^1$-smooth
with $f(0,0)=0$ and $g(0,0)=0$, and satisfy Lipschitz condition as
$$
\|f(x, y)-f(\tilde{x}, \tilde{y})\|_{H_s}\leq
K(\|x-\tilde{x}\|_{H_s}+\|y-\tilde{y}\|_{H_f}),
$$
$$
\|g(x, y)-g(\tilde{x}, \tilde{y})\|_{H_s}\leq
K(\|x-\tilde{x}\|_{H_s}+\|y-\tilde{y}\|_{H_f}),
$$
where  $K<\frac{1}{3}$. For example, $f(x^\varepsilon,
y^\varepsilon)=\frac{1}{4}\sin y^\varepsilon$ and $g(x^\varepsilon,
y^\varepsilon)=\frac{1}{4}\cos x^\varepsilon$.
\par
Then for the system (\ref{Equation-s-example
5.1})-(\ref{Equation-f-example 5.1}), taking
$\beta=\frac{-\gamma_s}{2}=\frac{1}{2}$,  as in Theorem 4.1 and
Theorem 5.1, the slow foliation and the critical foliation exist.
Their fibers are given as
$\mathcal{W}_{\frac{1}{2}}^\varepsilon((X_0, Y_0)^T, \omega)$ and
$\mathcal{W}_{\frac{1}{2}}^0((X_0, Y_0)^T, \omega)$, respectively.
Furthermore, the fiber $\mathcal{W}_{\frac{1}{2}}^\varepsilon((X_0,
Y_0)^T, \omega)$ of the slow foliation converges to the fiber
$\mathcal{W}_{\frac{1}{2}}^0((X_0, Y_0)^T, \omega)$ of  the critical
foliation in distribution, as the singular perturbation parameter
$\varepsilon$ tends to zero. In addition, the slow manifold
$\mathcal{M}^\varepsilon(\omega)$ as given by (\ref{manifold}) is
one fiber of the slow foliation, which parallels other fibers of the
slow foliation.}


\begin{thebibliography}{99}

\bibitem{Arnold}
L. Arnold, Random dynamical systems, Springer-Verlag, 1998.


\bibitem{BLZ}
P. Bates, K. Lu, and C. Zeng, Invariant foliations near normally
hyperbolic invariant manifolds for semiflows, \emph{Trans. Amer.
Math. Soc.}, \textbf{352}(10): 4641-4676, 2000.


\bibitem{BG-book}
N. Berglund, and B. Gentz, Noise-induced phenomena in slow-fast
dynamical systems: a sample-paths approach, Springer-Verlag, London,
2006.

\bibitem{BW}
D. Blomker, and W. Wang, Qualitative properties of local random
invariant manifolds for SPDE with quadratic nonlinearity, \emph{J.
Dyn. Differ. Equ.}, \textbf{22}: 677-695, 2010.


\bibitem{CDLS}
T. Caraballo, J. Duan, K. Lu, and B. Schmalfuss, Invariant manifolds
for random and stochastic partial differential equations, \emph{Adv.
Nonlinear Stud.}, \textbf{10}: 23-52, 2010.



\bibitem{CDZ}
G. Chen, J. Duan, and J. Zhang, Geometric shape of invariant
manifolds for a class of stochastic partial differential equations,
\emph{Journal of Mathematical Physics}, \textbf{52}: 072702, 2011.


\bibitem{CDZ2}
G. Chen, J. Duan, and J. Zhang, Approximating dynamics of a
singularly perturbed stochastic wave equation with a random
dynamical boundary condition, \emph{SIAM Journal on Mathematical
Analysis}, accepted, 2013.

\bibitem{CHT}
X. Chen, J. Hale, and B. Tan, Invariant foliations for $C^1$
semigroups in Banach spaces, \emph{J. Differential Equations},
\textbf{139}(2): 283-318, 1997.



\bibitem{CLL}
S.N. Chow, X.B. Lin, and K. Lu, Smooth invariant foliations in
infinite-dimensional spaces,  \emph{J. Diff. Eqs.} \textbf{94}:
266-291, 1991.

\bibitem{DZ}
G. Da Prato, and J. Zabczyk, Stochastic equations in infinite
dimensions, Cambridge University Press, Cambridge, England, 1992.


\bibitem{DLS-2003}
J. Duan, K. Lu, and B. Schmalfuss, Invariant manifolds for
stochastic partial differential equations, \emph{Ann. Probab.},
\textbf{31}(4): 2109-2135, 2003.

\bibitem{DLS-2004}
J. Duan, K. Lu, and B. Schmalfuss, Smooth stable and unstable
manifolds for stochastic evolutionary equations, \emph{J. Dyn.
Differ. Equ.}, \textbf{16}(4): 949-972, 2004.


\bibitem{Fenichel}
N. Fenichel, Geometric singular perturbation theory for ordinary
differential equations, \emph{J. Differential Equations},
\textbf{31}(1): 53-98, 1979.

\bibitem{Freidlin}
M.I. Freidlin, On stable oscillations and equilibriums induced by
small noise,\emph{ J. Statist. Phys.}, \textbf{103}: 283-300, 2001.


\bibitem{FLD}
H. Fu, X. Liu, and J. Duan, Slow manifolds for multi-time-scale
stochastic evolutionary systems, \emph{Commun. Math. Sci.},
\textbf{11}(1): 141-162, 2013.


\bibitem{Jones-book-1995}
C. K. R. T. Jones, Geometric singular perturbation theory, C.I.M.E. Lectures,
Montecatini Terme, June 1994, Lecture Notes in Mathematics 1609,
Springer-Verlag, Heidelberg, 1995.


\bibitem{KP-book}
Y. Kabanov, and S. Pergamenshchikov,  Two-scale stochastic systems:
asymptotic analysis and control, Springer-Verlag, Berlin, 2003.


\bibitem{LS08}
K. Lu, and B. Schmalfuss, Invariant foliations for stochastic
partial differential equations, \emph{Stochastics and Dynamics},
\textbf{8}(3): 505-518, 2008.


\bibitem{LS}
K. Lu, and B. Schmalfuss, Invariant manifolds for stochastic wave
equations, \emph{J. Differ. Equations}, \textbf{236}: 460-492, 2007.

\bibitem{MZZ}
S.A. Mohammed, T. Zhang, and H. Zhao, The stable manifold theorem
for semilinear stochastic evolution equations and stochastic partial
differential equations, \emph{ Mem. Am. Math. Soc.}, \textbf{196}:
1-105, 2008.



\bibitem{RDJ}
J. Ren, J. Duan, and C. K. R. T. Jones, Approximation of random slow
manifolds and settling of inertial particles under uncertainty,
http://arxiv.org/pdf/1212.4216v2.pdf, 2012.


\bibitem{SS}
B. Schmalfuss, and K. Schneider, Invariant manifolds for random
dynamical systems with slow and fast variables, \emph{Journal of
Dynamics and Differential Equations}, \textbf{20}(1): 133-163, 2008.


\bibitem{SKD}
X.Sun, X. Kan, and J. Duan, Approximation of invariant foliations
for stochastic dynamical systems, \emph{Stochastics and dynamics},
\textbf{12}(1): 1150011, 2012.


\bibitem{WRD-2012-JDE}
W. Wang,  A Roberts, and J. Duan, Large deviations and
approximations for slow-fast stochastic reaction-diffusion
equations, \emph{J. Differential Equations}, \textbf{253}:
3501-3522, 2012.

\bibitem{WR-2013-JMAA}
W. Wang, and A. Roberts, Slow manifold and averaging for slow-fast
stochastic differential systems, \emph{J. Math. Anal. Appl.},
\textbf{398}: 822-839, 2013.


\bibitem{WD}
E. Waymire, and J. Duan, Probability and partial differential
equations in modern applied mathematics, IMA Vol. 140, Springer, New
York, 2005.



\end{thebibliography}
\end{document}